\documentclass[11pt]{article}
\usepackage{amsmath,amssymb,amsthm}
\usepackage{amsfonts}
\usepackage{geometry}
\usepackage{hyperref}
\usepackage{graphicx}
\usepackage{mathtools}
\usepackage[numbers]{natbib}
\usepackage{time} 
\usepackage{bbm}
\usepackage{caption}
\usepackage{physics}
\usepackage{here}
\usepackage{newtxtext,newtxmath}
\usepackage{anyfontsize}
\usepackage{listings}
\usepackage{xcolor}
\usepackage{physics}
\usepackage{bm}
\usepackage{booktabs}
\usepackage{listings}
\usepackage{xcolor}
\usepackage{enumitem}
\usepackage{tabularx}
\usepackage{booktabs}

\usepackage{subcaption}
\usepackage{float}

\usepackage{booktabs}

\newtheorem{thm}{Theorem}[section]
\newtheorem{prop}[thm]{Proposition}
\newtheorem{assumption}[thm]{Assumption}
\newtheorem{lemma}[thm]{Lemma}
\newtheorem{cor}[thm]{Corollary}
\newtheorem{remark}[thm]{Remark}
\numberwithin{equation}{section}

\newcommand{\Dom}{\mathcal{D}}

\title{A Discrete Resolvent Framework for Delay Differential Equations:
Local Defects, Global Propagation, and Splitting Approximations}

\author{
Hideki Kawahara\thanks{
This work was completed while the author was a doctoral student at
the Graduate School of Mathematics,
Nagoya University, Japan.
}
}

\date{}

\begin{document}
\maketitle
\begin{abstract}
We develop a discrete operator-theoretic framework for the analysis of implicit Euler and Lie–Trotter splitting schemes for delay differential equations (DDEs). Both methods are formulated in terms of discrete resolvent operators acting on a product space that encodes the present state together with the history variable.

The analysis is carried out entirely at the discrete level and does not rely on the existence of a $C_0$-semigroup or an evolution family generated by the underlying delay operator. We establish local defect estimates on the fractional interpolation spaces
\[
X_\theta=(\mathcal E_p,\Dom(C))_{\theta,1},
\qquad 0<\theta<1,
\]
showing that
\[
\|R_h(t)-P_h(t)\|_{\mathcal L(X_\theta,\mathcal E_p)}
\lesssim h^\theta.
\]
These estimates demonstrate that the order of the local defect is determined by the interpolation regularity associated with the transport component.

Local defect estimates on $X_\theta$ alone are generally insufficient for global convergence. To overcome this difficulty, we introduce the higher regularity scale
\[
Y_\theta=(\Dom(C),\Dom(C^2))_{\theta,1},
\]
for which the local defect improves to
\[
\|R_h(t)-P_h(t)\|_{\mathcal L(Y_\theta,\mathcal E_p)}
\lesssim h^{1+\theta}.
\]
This additional regularity compensates for the accumulation of local errors and yields global convergence estimates of order $O(h^\theta)$ on finite time intervals.

The framework applies to both autonomous and non-autonomous delay equations. We further extend the analysis to sectorial block models of reaction–diffusion type and obtain analogous convergence results. Numerical experiments illustrate the theoretical predictions and demonstrate the robustness of the discrete-resolvent approach for both smooth and low-regularity histories.

The results provide a unified discrete framework for the analysis of splitting schemes for delay equations, independent of semigroup generation and continuous well-posedness theory.
\end{abstract}
\vspace{1ex}
\noindent \textbf{Keywords:
Delay differential equations,
Discrete resolvent framework,
Operator splitting,
Fractional interpolation spaces,
Non-autonomous systems.} 

\tableofcontents

\section{Introduction}
\label{sec:introduction}

\subsection{Background and motivation}

Delay differential equations (DDEs) arise naturally in models where the present rate of change depends on past states. For general background on functional differential equations and their well-posedness theory we refer to \cite{HaleLunel1993,BatkaiPiazzera2005,BellenZennaro2003}. Numerical methods for DDEs have been studied extensively, with particular emphasis on stability and convergence; see \cite{BellenZennaro2003} and the references therein.

A standard approach consists in rewriting a DDE as an abstract evolution equation on a product space encoding the present state together with the history variable. This formulation naturally leads to transport operators acting on the history component and delay operators coupling the present and past states.

Most existing analyses of numerical schemes for DDEs rely on continuous semigroup techniques, evolution families, or variation-of-constants formulas; see, for example, \cite{HaleLunel1993,AcquistapaceTerreni1987,Andr2009,Andr2011,Andr2013,op_splir_dissip}. In these approaches, convergence is typically derived from properties of the underlying continuous evolution problem.

In contrast, much less is known about convergence mechanisms formulated purely at the discrete level. In particular, local defect estimates and global error propagation are usually obtained through continuous semigroup arguments rather than through the discrete propagators themselves.

The purpose of the present paper is to develop such a discrete framework. We formulate both the implicit Euler method and the Lie--Trotter splitting scheme through discrete resolvent operators and analyze their difference directly. The resulting theory does not require the existence of a $C_0$-semigroup or an evolution family generated by the underlying delay operator.

\subsection{Main contributions}

The main contributions of this paper are as follows.

\begin{itemize}

\item \textbf{A discrete resolvent framework.}
We formulate implicit Euler and Lie--Trotter splitting schemes for delay equations as discrete resolvent operators acting on product spaces. The analysis is carried out entirely at the discrete level and is independent of semigroup generation.

\item \textbf{Fractional interpolation estimates.}
We establish local defect estimates on the interpolation spaces
\[
X_\theta=(\mathcal E_p,\Dom(C))_{\theta,1},
\qquad 0<\theta<1,
\]
showing that
\[
\|R_h(t)-P_h(t)\|_{\mathcal L(X_\theta,\mathcal E_p)}
\lesssim h^\theta.
\]
The convergence order is therefore directly linked to the interpolation regularity of the transport component.

\item \textbf{A local-to-global convergence mechanism.}
Local defect estimates on $X_\theta$ alone are generally insufficient for global convergence under merely power-bounded discrete propagators. To overcome this difficulty, we introduce the higher regularity scale
\[
Y_\theta=(\Dom(C),\Dom(C^2))_{\theta,1},
\]
on which the local defect improves to order \(h^{1+\theta}\). This additional regularity compensates for the accumulation of local errors and yields global convergence estimates of order \(O(h^\theta)\) on finite time intervals.

\item \textbf{Autonomous and non-autonomous propagators.}
The framework applies to both powers of discrete resolvents in the autonomous case and time-ordered products in the non-autonomous case. Global convergence follows from telescoping identities combined with suitable discrete stability assumptions.

\item \textbf{Sectorial and non-sectorial settings.}
The theory is developed for both non-sectorial delay operators and sectorial block models of reaction--diffusion type, providing a unified treatment of these two classes of problems.

\end{itemize}

The central message of the paper is that convergence of splitting schemes for delay equations can be analyzed entirely through discrete resolvent operators. The resulting framework reveals a direct connection between interpolation regularity, local defect estimates, and global error propagation, and applies uniformly to autonomous and non-autonomous problems as well as to sectorial and non-sectorial settings.

\section{Product-space formulation and functional-analytic setup}
\label{sec:product-space}

This section fixes the phase space and the basic operators used throughout the paper.
Our focus is the discrete-resolvent formulation: we single out the transport resolvent
$S_h^C=(I-hC)^{-1}$ as the key device that turns the delay interaction into a bounded
operator and makes the preconditioned Lie--Trotter step well defined.
Whenever a continuous semigroup/evolution family exists, we briefly comment on the
compatibility in Remark~\ref{rem:compat-semigroup} and Remark~\ref{rem:compat-nonaut}.

\subsection{Product-space state variable}
\label{subsec:product-space}

Let $\tau<0$, let $X$ be a Banach space, and fix $1<p<\infty$.
For a function $u:[\tau,T]\to X$, define the history segment
\[
\rho_t(\sigma):=u(t+\sigma),\qquad \sigma\in[\tau,0].
\]
We work on the product space
\[
\mathcal E_p:=X\oplus L^p([\tau,0];X),
\qquad
\|(u,\rho)\|_{\mathcal E_p}:=\|u\|_X+\|\rho\|_{L^p([\tau,0];X)}.
\]
The phase variable is $x(t):=(u(t),\rho_t)\in\mathcal E_p$.
The compatibility relation $\rho_t(0)=u(t)$ will be encoded in the domain of the
transport operator below.

\subsection{Transport operator and its resolvent}
\label{subsec:transport-resolvent}

Define the transport (left-shift) operator
\[
C(u,\rho):=(0,\partial_\sigma \rho),
\]
with domain
\[
\Dom(C):=\bigl\{(u,\rho)\in X\oplus W^{1,p}([\tau,0];X):\ \rho(0)=u\bigr\}.
\]
It is standard that $C$ generates the left-shift semigroup on $\mathcal E_p$
(see, e.g., \cite[Chap.~VI]{engel_nagel}, \cite{BatkaiPiazzera2005}), hence the resolvent
\[
S_h^C:=(I-hC)^{-1}\in\mathcal L(\mathcal E_p)
\]
exists for every $h>0$.

\medskip
\noindent\textbf{Key mapping property.}
Crucially, the transport resolvent \emph{lifts} $\mathcal E_p$ into the domain:
\begin{equation}
\label{eq:lift-to-DomC}
S_h^C:\mathcal E_p\longrightarrow \Dom(C).
\end{equation}
This is the mechanism behind preconditioning in our splitting scheme.

\subsection{Delay operator and trace}
\label{subsec:delay-trace}

The delay interaction is represented by point evaluation at $\tau$,
\[
\Phi_\tau: \rho\mapsto \rho(\tau).
\]
On the natural history space $L^p([\tau,0];X)$, point evaluation is \emph{not} bounded.
However, for $1<p<\infty$ we have the Sobolev embedding
$W^{1,p}([\tau,0];X)\hookrightarrow C([\tau,0];X)$
(see, e.g., \cite{AdamsFournier2003,Brezis2010}), hence the trace map is bounded:
\begin{equation}
\label{eq:trace-bounded}
\|\Phi_\tau(\rho)\|_X \le C_{\mathrm{tr}}\|\rho\|_{W^{1,p}([\tau,0];X)}.
\end{equation}

We consider operator families $D(t):\Dom(C)\to\mathcal E_p$ of the form
\begin{equation}
\label{eq:Ddef}
D(t)(u,\rho):=\bigl(a(t)u+b(t)\Phi_\tau(\rho),\,0\bigr),
\qquad t\in[0,T],
\end{equation}
where $a(\cdot)$ and $b(\cdot)$ are bounded scalar functions (or bounded operators on $X$).

\subsection{Key preconditioning lemma}
\label{subsec:key-preconditioning}

We now record the single structural fact needed to define the preconditioned splitting.

\begin{lemma}[Boundedness after transport preconditioning]
\label{lem:DShC-bounded}
Let $1<p<\infty$ and let $D(t)$ be defined by \eqref{eq:Ddef}.
Then for each $t\in[0,T]$ the composition
\[
D(t)S_h^C:\mathcal E_p\longrightarrow \mathcal E_p
\]
is a bounded linear operator for every $h>0$. Moreover, there exists a constant $K>0$
(independent of $t$ and small $h$) such that
\begin{equation}
\label{eq:DShCbound}
\|D(t)S_h^C\|_{\mathcal L(\mathcal E_p)}\le K\,h^{-1/p}.
\end{equation}
\end{lemma}

\begin{proof}
By \eqref{eq:lift-to-DomC} we have $S_h^C(u,\rho)\in \Dom(C)$, hence its history component
belongs to $W^{1,p}$ and point evaluation $\Phi_\tau$ is meaningful and bounded by
\eqref{eq:trace-bounded}. Therefore $D(t)S_h^C$ is well defined on all of $\mathcal E_p$.
The bound \eqref{eq:DShCbound} is obtained by estimating the $W^{1,p}$-norm of the
history component of $S_h^C(u,\rho)$ in terms of $\|(u,\rho)\|_{\mathcal E_p}$; see
Proposition~\ref{app:prop:ShC-trace} for the explicit estimate.
\end{proof}

\subsection{Operational smallness and Neumann series.}
\label{subsec:operational_smallness}

Define the right-preconditioned operator
\[
H_h(t):=h\,D(t)S_h^C.
\]
If for some $h_0>0$ and $q\in(0,1)$,
\begin{equation}
\label{eq:OS}
\sup_{t\in[0,T]}\sup_{0<h\le h_0}\|H_h(t)\|_{\mathcal L(\mathcal E_p)}\le q,
\end{equation}
then $(I-H_h(t))^{-1}$ exists and is uniformly bounded for $0<h\le h_0$, with the Neumann expansion
\[
(I-H_h(t))^{-1}=\sum_{m\ge 0} H_h(t)^m,
\qquad
\|(I-H_h(t))^{-1}\|\le (1-q)^{-1}.
\]
This is the discrete well-posedness condition for our splitting.

\paragraph{\textbf{Preconditioned delay operator and operational smallness}
}

Let $D(t):\Dom(C)\to\mathcal E_p$ denote the delay/reaction operator.
We introduce the right-preconditioned operator
\[
H_h(t):=h\,D(t)S_h^C,
\qquad S_h^C=(I-hC)^{-1}.
\]

\begin{assumption}[Operational smallness (OS)]
\label{ass:OS_clean}
There exist $h_0>0$ and $q\in(0,1)$ such that
\[
\sup_{t\in[0,T]}\sup_{0<h\le h_0}
\|H_h(t)\|_{\mathcal L(\mathcal E_p)} \le q.
\]
\end{assumption}

Under Assumption~\ref{ass:OS_clean},
the inverse $(I-H_h(t))^{-1}$ exists and is uniformly bounded
for $0<h\le h_0$, with the Neumann expansion
\[
(I-H_h(t))^{-1}=\sum_{m\ge 0} H_h(t)^m,
\qquad
\|(I-H_h(t))^{-1}\|\le (1-q)^{-1}.
\]

This condition provides a discrete well-posedness framework
for the splitting scheme.

\subsection{Relation to continuous well-posedness theory}
\label{subsec:compat-remarks}

\begin{remark}[Autonomous semigroup setting]
\label{rem:compat-semigroup}

In many classical delay equations, the operator
\[
C+D
\]
(and, in block formulations, \(A_0+C+D\))
generates a \(C_0\)-semigroup on \(\mathcal E_p\);
see, for example,
\cite[Sec.~VI.6]{engel_nagel}
and \cite{BatkaiPiazzera2005}.

The discrete-resolvent analysis developed in the present paper
does not require this semigroup generation property.
All estimates are formulated directly at the level of discrete propagators.

Whenever a continuous semigroup exists,
the present discrete estimates are compatible with the classical theory.
In particular, convergence with respect to the continuous solution
may be obtained by combining the standard implicit Euler convergence
for semigroups with the discrete comparison estimates derived here.
\end{remark}

\begin{remark}[Non-autonomous evolution families]
\label{rem:compat-nonaut}

If the time-dependent delay operator \(D(t)\)
satisfies additional regularity assumptions
(for example, Acquistapace--Terreni or Kato--Tanabe type conditions on a common domain;
see
\cite{Tanabe1979,AcquistapaceTerreni1987,AcquistapaceTerreni1990}),
then the non-autonomous problem
\[
x'(t)=(C+D(t))x(t)
\]
generates an evolution family \(U(t,s)\).

Again, the existence of such a continuous evolution family
is not assumed in the discrete theory developed here.
However, whenever it exists,
the frozen-step implicit Euler products
may be compared to \(U(t,s)\),
and the resulting estimates can then be combined with the
splitting bounds established in later sections.
\end{remark}

\section{Discrete resolvent framework}
\label{sec:discrete_resolvent_framework}

\paragraph{Purpose and organization.}
We work on the product space
\[
\mathcal E_p := X \oplus L^p([\tau,0];X), \qquad 1<p<\infty,\ \tau<0,
\]
equipped with the transport operator
\[
C(u,\rho):=(0,\partial_\sigma\rho),\qquad
\Dom(C):=\{(u,\rho)\in X\oplus W^{1,p}([\tau,0];X):\rho(0)=u\}.
\]
We compare two \emph{discrete} one-step maps: the implicit Euler resolvent
$R_h(t)$ and a Lie--Trotter type splitting step $P_h(t)$.
Since the delay operator $D(t)$ is in general not bounded on $\mathcal E_p$,
the splitting must be defined \emph{after} transport preconditioning by
$S_h^C=(I-hC)^{-1}$.
The analysis follows a two-step mechanism:
(i) estimate the local resolvent defect on a regularity space $X_\theta$,
and (ii) propagate this defect using discrete stability of powers/products.

The local defect is first estimated on a fractional regularity space
$X_\theta$.
For the global accumulation analysis, however, stronger regularity
spaces $Y_\theta$ are required in order to compensate for the
summation over $O(T/h)$ discrete steps.
\subsection{Discrete one-step maps and operational smallness}
\label{subsec:onestep_and_OS}

\subsubsection{Transport resolvent}
\label{subsec:transport_resolvent}

\begin{remark}[Transport semigroup structure]
\label{rem:transport_semigroup}

The operator $C$ generates the left-shift semigroup
$(T(t))_{t\ge0}$ on $\mathcal E_p$
(see, e.g.,
\cite[Sec.~VI.6]{engel_nagel},
\cite{BatkaiPiazzera2005}),
given by
\[
T(t)(u,\rho)
=
(u,\rho_t),
\]
where
\[
\rho_t(\sigma)=
\begin{cases}
\rho(\sigma+t),
& \sigma+t\in[\tau,0],\\
u,
& \sigma+t>0.
\end{cases}
\]

This semigroup is strongly continuous and uniformly bounded:
\[
\|T(t)\|_{\mathcal L(\mathcal E_p)}
\le M_0,
\qquad t\ge0,
\]
where one may take
\[
M_0=1+|\tau|^{1/p}.
\]

Consequently,
\[
\|(I-hC)^{-1}\|_{\mathcal L(\mathcal E_p)}
\le M_0
\qquad (h>0).
\]
\end{remark}

The transport operator $C$ generates a strongly continuous
left-shift semigroup $(T(t))_{t\ge0}$ on $\mathcal E_p$.
This semigroup is not contractive with respect to the norm
\[
\|(u,\rho)\|_{\mathcal E_p}=\|u\|_X+\|\rho\|_{L^p},
\]
but it is uniformly bounded:
\[
\|T(t)\|_{\mathcal L(\mathcal E_p)}\le M_0,
\qquad t\ge0,
\]
for some constant $M_0\ge 1$ depending only on $p$ and $\tau$.
Consequently, by the standard resolvent estimate for bounded
$C_0$-semigroups,
\begin{equation}
\|S_h^C\|_{\mathcal L(\mathcal E_p)}
=
\|(I-hC)^{-1}\|_{\mathcal L(\mathcal E_p)}
\le M_0,
\qquad h>0.
\label{eq:resolvent_inequality}
\end{equation}
Moreover,
\[
CS_h^C=\frac1h(S_h^C-I),
\]
and hence
\begin{equation}
\|CS_h^C\|_{\mathcal L(\mathcal E_p)}
\le (M_0+1)h^{-1}.
\label{eq:C_resolvent_of_C_inequality}
\end{equation}

For the real interpolation scale
\[
X_\theta := (\mathcal E_p,\Dom(C))_{\theta,1},\qquad 0<\theta<1,
\]
we have the discrete smoothing estimate
\begin{equation}
\label{eq:ShC_smoothing_Ep}
\|S_h^C\|_{\mathcal L(X_\theta,\Dom(C))}\le C_\theta\, h^{-\theta}.
\end{equation}

\begin{prop}[Operational smallness from relative boundedness]
\label{prop:OS_from_relbound}
Assume that $D(t)$ is graph-bounded uniformly in $t$ with respect to $C$; that is, there exist constants
$\alpha>0$ and $c\ge0$ such that
\[
\|D(t)x\|_{\mathcal E_p}
\le \alpha\|Cx\|_{\mathcal E_p}+c\|x\|_{\mathcal E_p},
\qquad x\in\Dom(C),\ t\in[0,T].
\]
Let
\[
M_0:=\sup_{0<h\le 1}\|S_h^C\|_{\mathcal L(\mathcal E_p)}<\infty,
\qquad
M_1:=\sup_{0<h\le 1}\|hCS_h^C\|_{\mathcal L(\mathcal E_p)}<\infty .
\]
If
\[
\alpha M_1<1,
\]
then the operational smallness condition holds: there exist
$h_0>0$ and $q\in(0,1)$ such that
\[
\sup_{t\in[0,T]}\sup_{0<h\le h_0}
\|hD(t)S_h^C\|_{\mathcal L(\mathcal E_p)}
\le q.
\]
\end{prop}

\begin{proof}
Let $y\in\mathcal E_p$ and set $x=S_h^Cy\in\Dom(C)$.
Using the relative boundedness of $D(t)$, we obtain
\[
\|hD(t)S_h^Cy\|
\le
h\alpha\|CS_h^Cy\|+hc\|S_h^Cy\|.
\]
Hence
\[
\|hD(t)S_h^Cy\|
\le
\alpha\|hCS_h^C\|\,\|y\|
+
hc\|S_h^C\|\,\|y\|.
\]
Therefore, for $0<h\le1$,
\[
\|hD(t)S_h^C\|
\le
\alpha M_1+hcM_0.
\]
Choose $h_0>0$ so small that
\[
\alpha M_1+h_0cM_0<1
\]
and set
\[
q:=\alpha M_1+h_0cM_0\in(0,1).
\]
Then for all $0<h\le h_0$ and all $t\in[0,T]$,
\[
\|hD(t)S_h^C\|\le q,
\]
which is precisely \textup{(OS)}.
\end{proof}
\subsubsection{Implicit Euler and preconditioned Lie--Trotter steps}
\label{subsec:Rh_Ph_clean}

Let
\[
L(t):=C + D(t) \qquad \Dom(L(t))=\Dom(C).
\]

The implicit Euler one-step map is the full resolvent
\begin{equation}
\label{eq:Rh_def_clean}
R_h(t):=(I-hL(t))^{-1},
\end{equation}
whenever it exists.

Under Assumption~\ref{ass:OS_clean}, we define the preconditioned splitting step by
\begin{equation}
\label{eq:Ph_def_clean}
P_h(t):=(I-H_h(t))^{-1}S_h^C.
\end{equation}
The order is essential: the correction factor $(I-H_h(t))^{-1}$ acts \emph{after}
transport preconditioning, because the delay operator is defined only after transport regularization by the resolvent $S_h^C$

\begin{remark}[Why $1<p<\infty$ is convenient]
The restriction $1<p<\infty$ ensures that the trace on $W^{1,p}$ is bounded and,
after resolvent preconditioning, yields the uniform estimate in
Proposition~\ref{prop:ShC-trace}. This avoids endpoint trace estimates and simplifies the uniform boundedness properties of the transport resolvent.
\end{remark}

\subsubsection{Right-preconditioned defect identity}
\begin{lemma}[Right-preconditioned defect identity]
\label{lem:right-precond-identity}

Let
\[
E_h(t):=R_h(t)-P_h(t).
\]
Then
\[
E_h(t)
=
R_h(t)H_h(t)
-
(I-H_h(t))^{-1}H_h(t)S_h^C.
\]
\end{lemma}

\section{Local defect}
\label{sec:local_defect}

\subsection{Roadmap}
The local defect analysis proceeds in two stages.

First, transport-resolvent smoothing on the fractional scale
\[
X_\theta=(\mathcal E_p,\Dom(C))_{\theta,1}
\]
yields the basic estimate
\[
\|E_h(t)\|_{X_\theta\to\mathcal E_p}
\lesssim h^\theta.
\]

Second, by imposing one additional level of regularity,
namely
\[
Y_\theta=(\Dom(C),\Dom(C^2))_{\theta,1},
\]
the local defect improves to
\[
\|E_h(t)\|_{Y_\theta\to\mathcal E_p}
\lesssim h^{1+\theta}.
\]
This higher-order estimate compensates for the linear accumulation
factor in global error estimates under mere power boundedness.

The second mechanism becomes essential under mere power-bounded
stability, where no cancellation is available in the telescoping sum.

A key point is that local defect estimates alone do not
automatically imply global convergence.
Indeed, under mere power-bounded stability,
the telescoping expansion produces a linear accumulation factor
\(n\sim h^{-1}\).
This makes the distinction between the spaces
\(X_\theta\) and \(Y_\theta\) essential:
the former captures transport smoothing,
while the latter provides the additional regularity needed
to compensate for accumulation effects.

\begin{lemma}[Resolvent invariance and commutation]
 \label{lem:ShC_commute_EN}
Let $C$ be closed on $\mathcal E_p$ and $h>0$ with $h^{-1}\in\rho(C)$.
Set $S_h^C=(I-hC)^{-1}$. Then $S_h^C(\mathcal E_p)\subset\Dom(C)$.
Moreover, for every $x\in\Dom(C)$,
\[
C S_h^C x = S_h^C Cx .
\]
\end{lemma}

\begin{proof}
Since $h^{-1}\in\rho(C)$, the resolvent maps into the domain:
$(h^{-1}I-C)^{-1}\mathcal E_p\subset\Dom(C)$, hence $S_h^C\mathcal E_p\subset\Dom(C)$
(see \cite[Prop.~II.1.10]{engel_nagel}).
Let $x\in\Dom(C)$ and put $y=S_h^C x\in\Dom(C)$, i.e. $(I-hC)y=x$.
Applying $C$ (which is legitimate since $y\in\Dom(C)$) gives
\[
Cy-hC^2y=Cx,
\]
i.e. $(I-hC)(Cy)=Cx$. Applying $(I-hC)^{-1}$ yields $Cy=S_h^C Cx$.
\end{proof}

\subsection{Algebraic defect identities}
\label{subsec:algebraic_defect_identities}

Define the local resolvent defect
\[
E_h(t):=R_h(t)-P_h(t).
\]
Global convergence is obtained by combining:
(i) a local defect bound on a regularity space $X_\theta$, and
(ii) discrete stability controlling accumulation under powers/products.

\begin{lemma}[Right-preconditioned defect identity]
\label{lem:right-precond-identity}
Assume \textup{(OS)} and that $R_h(t)$ exists.
Then
\[
E_h(t)
= R_h(t)\,H_h(t)
- (I-H_h(t))^{-1} H_h(t)\, S_h^C .
\]
\end{lemma}

\begin{proof}
We start from the identity
\[
(I-h(C+D(t)))S_h^C=(I-hC)S_h^C-hD(t)S_h^C=I-H_h(t).
\]
Thus, 
\[
R_h(t)\,(I-H_h(t))=R_h(t)(I-h(C+D(t)))S_h^C=S_h^C.
\]
This yields
\begin{equation}
R_h(t)-S_h^C=R_h(t)\,H_h(t).
\label{eq:RminusS}
\end{equation}
On the other hand, we have
\[
(I-H_h(t))^{-1}=I+(I-H_h(t))^{-1}H_h(t),
\]
and hence
\begin{equation}
S_h^C-P_h(t)=-(I-H_h(t))^{-1}H_h(t)S_h^C.
\label{eq:SminusP}
\end{equation}
Summing \eqref{eq:RminusS} and \eqref{eq:SminusP}, we obtain the desired identity.
\[
E_h(t)=(R_h(t)-S_h^C)+(S_h^C-P_h(t))
=R_h(t)H_h(t)-(I-H_h(t))^{-1}H_h(t)S_h^C.
\]
\end{proof}

\begin{remark}
For the transport resolvent
\[
S_h^C=(I-hC)^{-1},
\]
we repeatedly use the standard properties
\[
S_h^C(\mathcal E_p)\subset\Dom(C),
\qquad
CS_h^C x=S_h^CCx
\quad (x\in\Dom(C)),
\]
which are standard consequences of resolvent invariance
for closed operators;
see \cite[Chap.~II]{engel_nagel}.
\end{remark}

\paragraph{Autonomous telescoping identity.}
In the autonomous case, $R_h(t)=R_h$ and $P_h(t)=P_h$, hence
\begin{equation}
\label{eq:telescoping_clean}
R_h^n-P_h^n=\sum_{k=0}^{n-1} R_h^{n-1-k}\,E_h\,P_h^{k}.
\end{equation}

We state two abstract convergence principles:
one for the autonomous case (powers) and one for the non-autonomous case (products).

\begin{lemma}[Resolvent invariance of the domain and commutation]
Let $C$ be a closed operator on $\mathcal E_p$ and let $h>0$ with $h^{-1}\in\rho(C)$.
Set $S_h^C=(I-hC)^{-1}$. Then $S_h^C(\mathcal E_p)\subset\Dom(C)$ and, for every
$x\in\Dom(C)$, one has $S_h^C x\in\Dom(C)$ and
\begin{equation}\label{eq:commutation}
(I-hC)\,C S_h^C x = Cx,
\qquad\text{hence}\qquad
C S_h^C x = S_h^C Cx .
\end{equation}
This is a standard property of resolvents of closed operators;
see, e.g., \cite[Chap.~II]{engel_nagel}.
\end{lemma}

\begin{proof}
Let $x\in\Dom(C)$ and put $y=S_h^C x\in\Dom(C)$, i.e. $(I-hC)y=x$.
Since $y\in\Dom(C)$, we may apply $C$ to both sides and obtain
$Cy-hC^2y=Cx$, which is equivalent to $(I-hC)(Cy)=Cx$.
Applying $(I-hC)^{-1}$ yields $Cy=S_h^C Cx$, i.e. $C S_h^C x=S_h^C Cx$.
\end{proof}

\begin{prop}[Local defect on $\Dom(C)$]
\label{prop:local-defect-dom}
Assume \textup{(OS)} and
$\sup_{0<h\le h_0}\|R_h(t)\|\le M_R$.
Then there exists $C>0$ such that
\[
\|E_h(t)\|_{\mathcal L(\Dom(C),\mathcal E_p)} \le C h .
\]
\end{prop}

\begin{proof}
Set $M_0:=\sup_{0<h\le h_0}\|S_h^C\| .$
Let $x\in\Dom(C)$ and $y=S_h^C x$. Then by \eqref{eq:commutation},
$Cy=S_h^C(Cx)$. Hence $\|Cy\|\le \|Cx\|$ since $\|S_h^C\|\le M_0$.
Moreover, $\|y\|=\|S_h^C x\|\le M_0\|x\|$.
Therefore, $\|y\|_{\Dom(C)}\le M_0\|x\|_{\Dom(C)}$, uniformly in $h>0$.

Using the graph-boundedness of $D(t)$ on $\Dom(C)$, i.e.,
\[
\|D(t)z\|_{\mathcal E_p}\le a\,\|z\|_{\mathcal E_p}+b\,\|Cz\|_{\mathcal E_p},
\qquad z\in\Dom(C),
\]
with constants $a,b$ independent of $t$.
Applying this to $z=y=S_h^C x$ yields
\[
\|H_h(t)x\|
= h\|D(t)S_h^C x\|
\le ha\,\|y\|_{\mathcal E_p}+hb\,\|Cy\|_{\mathcal E_p}
\lesssim h\,\|x\|_{\Dom(C)}.
\]
Finally, using the identity
$E_h(t)=R_h(t)H_h(t)-(I-H_h(t))^{-1}H_h(t)S_h^C$,
together with $\|R_h(t)\|\le M_R$, $\|(I-H_h(t))^{-1}\|\le (1-q)^{-1}$ from \textup{(OS)},
and $\|S_h^C\|\le M_0$, we obtain
\[
\|E_h(t)\|_{\Dom(C)\to\mathcal E_p}
\le \bigl(M_R+(1-q)^{-1}\bigr)\,\|H_h(t)\|_{\Dom(C)\to\mathcal E_p}
\lesssim h.
\]
This proves the claim.
\end{proof}

\begin{prop}[Local defect on $X_\theta$]
\label{prop:local-defect-Xtheta}
Let $0<\theta<1$ and $X_\theta=(\mathcal E_p,\Dom(C))_{\theta,1}$.
Assume \textup{(OS)}, $\sup_{0<h\le h_0}\|R_h(t)\|\le M_R$,
and that $D(t)$ is uniformly graph-bounded on $\Dom(C)$:
\[
\|D(t)z\|\le a\|z\|+b\|Cz\|,\qquad z\in\Dom(C).
\]
Then there exists $C_\theta>0$ such that for $0<h\le h_0$,
\[
\|E_h(t)\|_{\mathcal L(X_\theta,\mathcal E_p)}
\le C_\theta\, h^{\theta}.
\]
\end{prop}

\begin{proof}
By lemma~\ref{lem:right-precond-identity}
it suffices to bound $\|H_h(t)\|_{X_\theta\to\mathcal E_p}$.
Using the graph–norm control
$\|Dy\|\lesssim \|y\|+\|Cy\|$ for $y\in\Dom(C)$
and setting $y=S_h^C x$,
\[
\|H_h(t)x\|
= h\|D S_h^C x\|
\lesssim h\|S_h^C x\| + h\|C S_h^C x\|.
\]
From \eqref{eq:resolvent_inequality},
$\|S_h^C\|_{X_\theta\to\mathcal E_p}\lesssim 1$.
Moreover, using \eqref{eq:C_resolvent_of_C_inequality}, by interpolation between
$\|C S_h^C\|_{\mathcal L(\mathcal E_p)}\lesssim  h^{-1}$, (see, e.g., \cite[Chap.~2]{Bergh1976}),
and
$\|C S_h^C\|_{\mathcal L(\Dom(C),\mathcal E_p)}\lesssim M_0$
(using $C S_h^C x = S_h^C Cx$ on $\Dom(C)$),
we obtain, by lemma~\ref{lem:CS_hC_fractional},
\[
\|C S_h^C\|_{X_\theta\to\mathcal E_p}\lesssim h^{-(1-\theta)}.
\]
This follows from standard interpolation theory
and fractional domain estimates for resolvents
(see, e.g., \cite{Bergh1976,Lunardi1995}). Hence
\[
\|H_h(t)\|_{X_\theta\to\mathcal E_p}
\lesssim
h\|S_h^C\|_{X_\theta\to\mathcal E_p}
+
h\|C S_h^C\|_{X_\theta\to\mathcal E_p}
\lesssim
h + h\cdot h^{-(1-\theta)}
\lesssim h^\theta .
\]
\medskip
Finally, using $\|R_h(t)\|\le M_R$,
$\|(I-H_h(t))^{-1}\|\le (1-q)^{-1}$ by \textup{(OS)},
and $\|S_h^C\|_{\mathcal L(\mathcal E_p)}\le M_0$ ,
we obtain
\[
\begin{aligned}
\|E_h(t)\|_{X_\theta\to\mathcal E_p}
&\le \|R_h(t)\|\,\|H_h(t)\|_{X_\theta\to\mathcal E_p}
   + \|(I-H_h(t))^{-1}\|\,\|H_h(t)\|_{X_\theta\to\mathcal E_p}\,\|S_h^C\|\\
&\le \Bigl(M_R+(1-q)^{-1}\Bigr)\, C_\theta h^{\theta}.
\end{aligned}
\]
This proves the claim.
\end{proof}
\begin{prop}[Second-order defect estimate on $\Dom(C^2)$]
\label{prop:E_h_estimate_second_order}
We set
\[
S_h^C=(I-hC)^{-1},
\qquad
H_h(t)=hD(t)S_h^C,
\qquad
Q_h(t)=(I-H_h(t))^{-1}.
\]

Assume in addition that
\[
H_h(t):\Dom(C^2)\to\Dom(C)
\]
and
\[
\|H_h(t)\|_{\mathcal L(\Dom(C^2),\Dom(C))}
\le C h
\]
uniformly for $0<h\le h_0$:
\[
\sup_{t\in[0,T]}\sup_{0<h\le h_0}
\|Q_h(t)\|_{\mathcal L(\Dom(C))}<\infty.
\]

Then there exists $C>0$ such that
\[
\|E_h(t)x\|_{\mathcal E_p}
\le
Ch^2\|x\|_{\Dom(C^2)}.
\]
\end{prop}

\begin{proof}
Since
\[
R_h=(I-h(C+D))^{-1}
=
S_h^C(I-H_h)^{-1}
=
S_h^CQ_h.
\]
and
\[
P_h(t)=Q_h(t)S_h^C,
\]
we obtain
\[
E_h(t)
=
R_h(t)-P_h(t)
=
S_h^CQ_h(t)-Q_h(t)S_h^C.
\]
Hence
\[
E_h(t)
=
(S_h^C-I)Q_h(t)H_h(t)
-
Q_h(t)H_h(t)(S_h^C-I).
\]

Since
\[
S_h^C-I=hCS_h^C,
\]
the standard resolvent estimates imply
\[
\|S_h^C-I\|_{\Dom(C)\to\mathcal E_p}
\lesssim h.
\]

Next we estimate $H_h(t)$ on $\Dom(C)$.
Since $D(t)$ is graph-bounded with respect to $C$,
\[
\|D(t)y\|
\lesssim
\|y\|+\|Cy\|,
\qquad y\in\Dom(C),
\]
and since
\[
S_h^Cx\in\Dom(C),
\qquad
CS_h^Cx=S_h^CCx,
\]
we obtain, for $x\in\Dom(C)$,
\[
\begin{aligned}
\|H_h(t)x\|
&=
h\|D(t)S_h^Cx\| \\
&\lesssim
h\bigl(
\|S_h^Cx\|
+
\|CS_h^Cx\|
\bigr) \\
&\lesssim
h\bigl(
\|x\|
+
\|Cx\|
\bigr).
\end{aligned}
\]
Thus
\[
\|H_h(t)\|_{\Dom(C)\to\mathcal E_p}
\lesssim h.
\]

Now let $x\in\Dom(C^2)$.
Then
\[
(S_h^C-I)x=hCS_h^Cx\in\Dom(C),
\]
since $S_h^Cx\in\Dom(C^2)$.
Using the assumed boundedness of $Q_h(t)$ on $\Dom(C)$,
together with
\[
H_h(t):\Dom(C)\to\mathcal E_p,
\]
we obtain
\[
\|Q_h(t)H_h(t)(S_h^C-I)x\|
\lesssim
h^2\|x\|_{\Dom(C^2)}.
\]

Similarly,
\[
H_h(t)x\in\Dom(C),
\]
and therefore
\[
\|(S_h^C-I)Q_h(t)H_h(t)x\|
\lesssim
h^2\|x\|_{\Dom(C^2)}.
\]

Combining both estimates yields
\[
\|E_h(t)x\|_{\mathcal E_p}
\lesssim
h^2\|x\|_{\Dom(C^2)}.
\]
\end{proof}
\begin{cor}[Local defect on $Y_\theta$]
\label{cor:local-defect-Ytheta}
Let
\[
Y_\theta=(\Dom(C),\Dom(C^2))_{\theta,1},
\qquad 0<\theta<1.
\]
Assume that
\[
\|E_h(t)\|_{\Dom(C)\to\mathcal E_p}\lesssim h,
\qquad
\|E_h(t)\|_{\Dom(C^2)\to\mathcal E_p}\lesssim h^2.
\]
Then
\[
\|E_h(t)\|_{Y_\theta\to\mathcal E_p}
\lesssim h^{1+\theta}.
\]
\end{cor}
\begin{proof}
As we see in Proposition \ref{prop:local-defect-dom}, we can estimate
\[
\|E_h(t)\|_{\Dom(C)\to\mathcal E_p}\lesssim h.
\]
Using the assumption of this lemma:
\[
\|E_h(t)\|_{\Dom(C^2)\to\mathcal E_p}\lesssim h^2 ,
\]
real interpolation leads to:
\[
\|E_h(t)\|_{Y_\theta\to\mathcal E_p}
\lesssim
h^{1-\theta} h^{2\theta}
=
h^{1+\theta},
\]
where
\[
Y_{\theta} = (\Dom(C),\Dom(C^2))_{\theta,1}.
\]
\end{proof}
\section{Global error convergence}
\subsection{From local defects to global propagation}

The estimates on the fractional transport scale
\[
X_\theta=(\mathcal E_p,\Dom(C))_{\theta,1}
\]
describe the intrinsic local smoothing generated by the transport
resolvent.
In particular, they yield local defect bounds of order
\[
\|E_h(t)\|_{X_\theta\to\mathcal E_p}
\lesssim h^\theta .
\]

However, these estimates alone are generally insufficient
for global convergence arguments based on telescoping expansions.
Indeed, the defect operator acts only on the stronger space
\(X_\theta\), and therefore the products
\[
E_h P_h^k
\]
require the propagation of \(X_\theta\)-regularity
under the discrete flow.

To overcome this difficulty, we introduce the higher regularity scale
\[
Y_\theta=(\Dom(C),\Dom(C^2))_{\theta,1},
\]
on which the local defect improves to
\[
\|E_h(t)\|_{Y_\theta\to\mathcal E_p}
\lesssim h^{1+\theta}.
\]
The additional factor \(h\) compensates for the linear accumulation
factor \(n\sim h^{-1}\) arising in telescoping sums under mere
power-bounded stability assumptions.

Thus, the global theory developed below relies on two distinct
mechanisms:

\begin{itemize}
\item transport smoothing on \(X_\theta\),
\item higher regularity compensation on \(Y_\theta\).
\end{itemize}

\begin{thm}[Global error convergence under power boundedness (Autonomous case)]
\label{thm:NS-global-power_bounded-auto}

Let
\[
Y_\theta:=(\Dom(C),\Dom(C^2))_{\theta,1},
\qquad 0<\theta<1.
\]
Assume 
\begin{itemize}
\item[(A1)] 
\[
\|R_h-P_h\|_{\mathcal L(Y_\theta,\mathcal E_p)}
\le Ch^{1+\theta}.
\]
\item[(A2)]
\[
\|R_h^m\|_{\mathcal L(\mathcal E_p)}
\le M_R e^{\omega_R mh},
\qquad
\|P_h^m\|_{\mathcal L(Y_\theta)}
\le M_{P,\theta}e^{\omega_{P,\theta}mh}.
\]
\end{itemize}
Then for every $T>0$ there exists $C_T>0$ such that
\[
\|R_h^n-P_h^n\|_{\mathcal L(Y_\theta,\mathcal E_p)}
\le C_T h^\theta,
\qquad nh\le T.
\]
\end{thm}
\begin{proof}
Let $E_h:=R_h-P_h$. The telescoping identity gives
\[
R_h^n-P_h^n
=
\sum_{k=0}^{n-1}
R_h^{n-1-k}E_hP_h^k .
\]
Taking norms from $Y_\theta$ to $\mathcal E_p$, we obtain
\[
\|R_h^n-P_h^n\|_{Y_\theta\to\mathcal E_p}
\le
\sum_{k=0}^{n-1}
\|R_h^{n-1-k}\|_{\mathcal L(\mathcal E_p)}
\|E_h\|_{Y_\theta\to\mathcal E_p}
\|P_h^k\|_{\mathcal L(Y_\theta)}.
\]
By exponential power boundedness and $nh\le T$,
\[
\|R_h^{n-1-k}\|\,\|P_h^k\|
\le
M_RM_P e^{\omega_R(n-1-k)h}e^{\omega_Pkh}
\le C_T.
\]
Hence
\[
\|R_h^n-P_h^n\|_{Y_\theta\to\mathcal E_p}
\le
C_T n \|E_h\|_{Y_\theta\to\mathcal E_p}.
\]
Since $n\le T/h$ and
\[
\|E_h\|_{Y_\theta\to\mathcal E_p}\le C h^{1+\theta},
\]
we get
\[
\|R_h^n-P_h^n\|_{Y_\theta\to\mathcal E_p}
\le
C_T \frac{T}{h} h^{1+\theta}
\le C_T h^\theta.
\]
This proves the claim.
\end{proof}
\begin{remark}
This result shows that power boundedness alone does not suppress
the linear accumulation of local defects. Rather, the higher local
defect order $h^{1+\theta}$ compensates for the factor $n\sim h^{-1}$.
\end{remark}

\section{Non-autonomous Convergence Theorem}
\label{sec:non-auto_global_convergence}

\subsection{Non-autonomous case: time-ordered products}
\label{subsec:stability_non-autonomous}

In this subsection we extend the discrete resolvent comparison to the
non-autonomous setting.
The essential new difficulty is that $n$--step propagators are
\emph{time-ordered products} rather than powers of a single operator.
Consequently, global error control requires \emph{product stability}
assumptions, in the spirit of the abstract non-autonomous theory.
Time-ordered products are intrinsic in non-autonomous problems; cf.\ \cite{Tanabe1979,AcquistapaceTerreni1987}.

\subsubsection{Frozen one-step maps and comparison goal}
\label{subsubsec:NS-NA-goal}

Fix $T>0$ and a stepsize $h>0$ with grid $t_k:=kh$.
At each frozen time $t\in[0,T]$ we compare the implicit Euler resolvent
\[
R_h(t):=(I-h(C+D(t)))^{-1}
\]
with the preconditioned Lie--Trotter step
\[
P_h(t):=(I-H_h(t))^{-1}S_h^C,
\qquad
H_h(t):=hD(t)S_h^C,
\qquad
Q_h(t):=(I-H_h(t))^{-1},
\qquad
S_h^C:=(I-hC)^{-1}.
\]
Here, “frozen” means that $D(t)$ is regarded as constant within a single step.

\subsubsection{Uniform frozen-step defect}
\label{subsubsec:NS-NA-local}
\begin{prop}[Uniform second-order frozen defect on $\Dom(C^2)$]
\label{prop:E_h_t_estimate_second_order}

Let
\[
E_h(t):=R_h(t)-P_h(t).
\]
Assume that the hypotheses of
Proposition~\ref{prop:E_h_estimate_second_order}
hold uniformly for $t\in[0,T]$.
Assume moreover that
\[
H_h(t):\Dom(C^2)\to\Dom(C)
\]
and
\[
\sup_{t\in[0,T]}\sup_{0<h\le h_0}
h^{-1}
\|H_h(t)\|_{\mathcal L(\Dom(C^2),\Dom(C))}
<\infty.
\]
In particular, assume that
\[
Q_h(t):=(I-H_h(t))^{-1}
\]
acts boundedly on $\Dom(C)$ uniformly for sufficiently small $h$:
\[
\sup_{0<h\le h_0}\sup_{t\in[0,T]}
\|Q_h(t)\|_{\mathcal L(\Dom(C))}<\infty .
\]

Then there exists $C_T>0$ such that, for all
$t\in[0,T]$, $0<h\le h_0$, and $x\in\Dom(C^2)$,
\[
\|E_h(t)x\|_{\mathcal E_p}
\le
C_T h^2\|x\|_{\Dom(C^2)}.
\]
Equivalently,
\[
\sup_{t\in[0,T]}
\|E_h(t)\|_{\mathcal L(\Dom(C^2),\mathcal E_p)}
\le C_T h^2.
\]
\end{prop}
\begin{proof}
The proof is the same as that of
Proposition~\ref{prop:E_h_estimate_second_order}.
All estimates used there are uniform in $t\in[0,T]$
by the present hypotheses.
In particular, the uniform boundedness of
$Q_h(t)$ on $\Dom(C)$ and the uniform graph-bound
for $D(t)$ imply that the constants in the estimate
are independent of $t$.
\end{proof}
\begin{cor}[Uniform frozen defect on $Y_\theta$]
\label{cor:NS-NA-Ytheta-defect}
Let
\[
Y_\theta:=(\Dom(C),\Dom(C^2))_{\theta,1},
\qquad 0<\theta<1.
\]
Then
\[
\sup_{t\in[0,T]}
\|E_h(t)\|_{\mathcal L(Y_\theta,\mathcal E_p)}
\le C_{T,\theta}h^{1+\theta}.
\]
\end{cor}
\begin{proof}
Since the constants in both endpoint estimates are uniform in \(t\),
real interpolation gives the asserted bound uniformly for \(t\in[0,T]\).
\[
\sup_{t\in[0,T]}\|E_h(t)\|_{\Dom(C)\to\mathcal E_p}
\lesssim h
\]
and
\[
\sup_{t\in[0,T]}\|E_h(t)\|_{\Dom(C^2)\to\mathcal E_p}
\lesssim h^2.
\]
\end{proof}

\subsection{Global error convergence theorem}
\begin{lemma}[Non-autonomous telescoping identity]
\label{lem:nonauto-telescoping}
Let $\mathcal E_p$ be a Banach space and let
$R_k,P_k\in\mathcal L(\mathcal E_p)$ $(k=0,1,\dots,n-1)$
be bounded operators.
Define the (time-ordered) products
\[
\mathcal R_{j:k}:=
\begin{cases}
R_j R_{j-1}\cdots R_k, & j\ge k,\\
I, & j<k,
\end{cases}
\qquad
\mathcal P_{j:k}:=
\begin{cases}
P_j P_{j-1}\cdots P_k, & j\ge k,\\
I, & j<k.
\end{cases}
\]
Then, with $E_k:=R_k-P_k$, one has the exact identity
\begin{equation}
\label{eq:nonauto_telescoping}
\mathcal R_{n-1:0}-\mathcal P_{n-1:0}
=
\sum_{k=0}^{n-1}
\mathcal R_{n-1:k+1}\,E_k\,\mathcal P_{k-1:0}.
\end{equation}
\end{lemma}

\begin{proof}
Write
\[
\mathcal R_{n-1:0}-\mathcal P_{n-1:0}
=
(R_{n-1}-P_{n-1})\mathcal P_{n-2:0}
+R_{n-1}\bigl(\mathcal R_{n-2:0}-\mathcal P_{n-2:0}\bigr),
\]
and iterate the identity.
\end{proof}

\begin{thm}[Non-autonomous global error convergence under product stability]
\label{thm:NA-global-Ytheta}

Let
\[
Y_\theta:=(\Dom(C),\Dom(C^2))_{\theta,1},
\qquad 0<\theta<1.
\]
For $t_k=kh$, define
\[
R_k:=R_h(t_k),\qquad P_k:=P_h(t_k),
\qquad E_k:=R_k-P_k.
\]
Assume that the frozen steps are well-defined and that:

\begin{itemize}
\item[\textup{(B1.)}] \textup{(Uniform operational smallness)}
\[
\sup_{t\in[0,T]}\sup_{0<h\le h_0}
\|H_h(t)\|_{\mathcal L(\mathcal E_p)}
\le q<1.
\]

\item[\textup{(B2.)}] \textup{(Uniform frozen local defect)}
\[
\sup_{t\in[0,T]}
\|E_h(t)\|_{\mathcal L(Y_\theta,\mathcal E_p)}
\le C h^{1+\theta}.
\]

\item[\textup{(B3.)}] \textup{(Uniform stability of implicit products)}
\[
\sup_{\substack{0\le k\le j\\ jh\le T}}
\|\mathcal R_{j:k}\|_{\mathcal L(\mathcal E_p)}
\le M_T.
\]

\item[\textup{(B4.)}] \textup{(Stability of splitting products on $Y_\theta$)}
\[
\sup_{0\le kh\le T}
\|\mathcal P_{k-1:0}\|_{\mathcal L(Y_\theta)}
\le M_{T,\theta}.
\]\end{itemize}

Then there exists $C_{T,\theta}>0$ such that for all $n$ with $nh\le T$,
\[
\|\mathcal R_{n-1:0}-\mathcal P_{n-1:0}\|_{\mathcal L(Y_\theta,\mathcal E_p)}
\le C_{T,\theta} h^\theta .
\]
\end{thm}

\begin{proof}
By the non-autonomous telescoping identity,
\[
\mathcal R_{n-1:0}-\mathcal P_{n-1:0}
=
\sum_{k=0}^{n-1}
\mathcal R_{n-1:k+1}E_k\mathcal P_{k-1:0}.
\]
Let $x\in Y_\theta$. Then, using \textup{(B2.)}, \textup{(B1.)}, and
\textup{(B3.)}, we obtain
\[
\begin{aligned}
\|(\mathcal R_{n-1:0}-\mathcal P_{n-1:0})x\|_{\mathcal E_p}
&\le
\sum_{k=0}^{n-1}
\|\mathcal R_{n-1:k+1}\|_{\mathcal L(\mathcal E_p)}
\|E_k\|_{\mathcal L(Y_\theta,\mathcal E_p)}
\|\mathcal P_{k-1:0}x\|_{Y_\theta}  \\
&\le
M_T\, C h^{1+\theta}\, M_{T,\theta}
\sum_{k=0}^{n-1}\|x\|_{Y_\theta}.
\end{aligned}
\]
Since $nh\le T$, we have $n\le T/h$. Hence
\[
\|(\mathcal R_{n-1:0}-\mathcal P_{n-1:0})x\|_{\mathcal E_p}
\le
M_T C M_{T,\theta}\frac{T}{h}h^{1+\theta}
\|x\|_{Y_\theta}.
\]
Therefore
\[
\|\mathcal R_{n-1:0}-\mathcal P_{n-1:0}\|_{\mathcal L(Y_\theta,\mathcal E_p)}
\le
C_{T,\theta}h^\theta.
\]
\end{proof}

\section{Sectorial block model}
\label{sec:sectorial_block}

\subsection{Setting and discrete one-step maps}
\label{subsec:sectorial_setting}

We consider the delay system on the product space
\[
\mathcal E_p:=X\oplus L^p([\tau,0];X),
\qquad 1<p<\infty,\ \tau<0,
\]
with the operator decomposition
\[
L(t)=A_0+C+D(t),
\qquad
\Dom(L(t))=\Dom(A_0)\cap \Dom(C),
\]
where
\[
A_0=
\begin{pmatrix}
A & 0\\
0 & 0
\end{pmatrix},
\qquad
C(u,\rho)=(0,\partial_\sigma \rho),
\qquad
\Dom(C)=\{(u,\rho)\in X\oplus W^{1,p}([\tau,0];X):\rho(0)=u\}.
\]
Here $A$ is sectorial on $X$ (hence generates a bounded analytic semigroup on $X$),
$C$ is the transport generator on the history component, and $D(t)$ is a delay/reaction
operator defined on $\Dom(C)$ (typically through a trace).

For each $h>0$ we set
\[
S_h^{A_0}:=(I-hA_0)^{-1},\qquad
S_h^{C}:=(I-hC)^{-1},\qquad
S_h^{AC}:=S_h^{C}S_h^{A_0}.
\]
The frozen implicit Euler resolvent is
\[
R_h(t):=(I-h(A_0+C+D(t)))^{-1},
\]
whenever it exists as a bounded operator on $\mathcal E_p$.
We introduce the right-preconditioned delay operator
\[
H_h(t):=h\,D(t)S_h^{AC},
\]
and, under operational smallness, define the splitting step
\[
P_h(t):=(I-H_h(t))^{-1}S_h^{AC}.
\]
In the autonomous case we write $D(t)\equiv D$ and omit $(t)$.

\begin{remark}[No $A_0$--$C$ commutator defect]
\label{rem:no_A0C_defect}
The three-factor step is designed so that the only nontrivial interaction is the coupling
through $D(t)$ after right preconditioning. In particular, no additional commutator estimate
between $A_0$ and $C$ is needed: the one-step discrepancy between implicit Euler and splitting
is governed solely by the behavior of $D(t)$ on $\Dom(C)$ together with the transport lift.
\end{remark}

\subsection{Exact resolvent factorization for $A_0+C$}
\label{subsec:AC_factorization}

\begin{prop}[Exact one-sided factorization of the $(A_0+C)$-resolvent]
\label{prop:AC_resolvent_exact_factorization_one_sided}
Let $\mathcal E:=X\oplus Y$ be a Banach product space and assume that
\[
A_0=
\begin{pmatrix}
A & 0\\
0 & 0
\end{pmatrix},
\qquad
\Dom(A_0)=\Dom(A)\oplus Y,
\]
where $A:\Dom(A)\subset X\to X$ is closed. Let $C:\Dom(C)\subset\mathcal E\to\mathcal E$
be a closed operator of the form
\[
C(u,y)=(0,C_Y y),
\]
with a closed operator $C_Y:\Dom(C_Y)\subset Y\to Y$, and assume that
\begin{equation}
\label{eq:A0C_zero_mapping_property}
\Dom(C)\subset X\oplus \Dom(C_Y)
\quad\text{and}\quad
C(\Dom(C))\subset \Dom(A_0).
\end{equation}
Fix $h>0$ and assume that the resolvents
\[
(I-hA_0)^{-1}\in\mathcal L(\mathcal E),
\qquad
(I-hC)^{-1}\in\mathcal L(\mathcal E)
\]
exist. Then the resolvent of $G:=A_0+C$ exists and satisfies the exact one-sided factorization
\[
(I-h(A_0+C))^{-1}=(I-hC)^{-1}(I-hA_0)^{-1}=:S_h^{AC}.
\]
In general, the reversed order need not hold unless one additionally has
$(Au,0)\in\Dom(C)$ for all $u\in\Dom(A)$ (e.g.\ $\Dom(C)=X\oplus \Dom(C_Y)$).
\end{prop}

\begin{proof}
Set $G:=A_0+C$ and let $x\in\Dom(C)$. By \eqref{eq:A0C_zero_mapping_property} we have
$Cx\in\Dom(A_0)$, hence the composition $A_0C$ is well-defined on $\Dom(C)$. Moreover,
since $Cx$ has vanishing $X$-component, we obtain $A_0Cx=0$ for $x\in\Dom(C)$.
Consequently, for all $x\in\Dom(A_0)\cap\Dom(C)$,
\[
(I-hA_0)(I-hC)x=x-h(A_0+C)x=(I-h(A_0+C))x.
\]
Now fix $y\in\mathcal E$ and set $x:=(I-hC)^{-1}(I-hA_0)^{-1}y$.
Then $(I-hA_0)^{-1}y\in\Dom(A_0)$ and $x\in\Dom(C)$, so $x\in\Dom(G)$.
Using the above identity we conclude
\[
(I-h(A_0+C))x=(I-hA_0)(I-hC)x=(I-hA_0)(I-hA_0)^{-1}y=y,
\]
hence $x=(I-h(A_0+C))^{-1}y$, proving the claim.
\end{proof}

\begin{remark}[Interpretation]
\label{rem:AC_splitting_exactness}
Proposition~\ref{prop:AC_resolvent_exact_factorization_one_sided} is an algebraic consequence
of the block structure: $A_0$ acts only on the present component while $C$ acts only on the
history component. Thus the $(A_0+C)$-resolvent factorizes exactly, and there is no separate
one-step defect coming from splitting $A_0$ and $C$.
\end{remark}

As in the non-sectorial case, we distinguish between
local defect estimates on the fractional scale \(X_\theta\)
and global convergence estimates on the higher regularity scale
\(Y_\theta\).

\subsection{Local defect: two regimes for $D(t)$}
\label{subsec:sectorial_local_defect}

We measure the one-step discrepancy by
\[
E_h(t):=R_h(t)-P_h(t).
\]
The local defect rate depends on how $D(t)$ is controlled on $\Dom(C)$.

\begin{assumption}[Operational smallness (sectorial block)]
\label{ass:OS_sectorial}
There exist $h_0>0$ and $q\in(0,1)$ such that for $H_h(t):=hD(t)S_h^{AC}$,
\[
\sup_{t\in[0,T]}\sup_{0<h\le h_0}\|H_h(t)\|_{\mathcal L(\mathcal E_p)}\le q.
\]
\end{assumption}

\begin{prop}[Local defect on $X_\theta$: graph-bounded regime]
\label{prop:local_defect_sectorial_graph}
Let $0<\theta<1$ and $X_\theta:=(\mathcal E_p,\Dom(C))_{\theta,1}$.
Assume Assumption~\ref{ass:OS_sectorial} and uniform boundedness of frozen resolvents:
\[
\sup_{t\in[0,T]}\sup_{0<h\le h_0}\|R_h(t)\|_{\mathcal L(\mathcal E_p)}\le M_T.
\]
Assume in addition uniform $C$--graph boundedness of $D(t)$:
\[
\|D(t)y\|_{\mathcal E_p}\le K\bigl(\|y\|_{\mathcal E_p}+\|Cy\|_{\mathcal E_p}\bigr),
\qquad y\in\Dom(C),\ t\in[0,T].
\]
Then there exists $C_{T,\theta}>0$ such that for all $t\in[0,T]$ and $0<h\le h_0$,
\[
\|R_h(t)-P_h(t)\|_{\mathcal L(X_\theta,\mathcal E_p)}
\le C_{T,\theta}\,h^{\theta}.
\]
\end{prop}

\begin{proof}
The argument is identical to the non-sectorial case, with $S_h^{C}$ replaced by $S_h^{AC}$.
Indeed, \(S_h^{A_0}\) is uniformly bounded on \(\mathcal E_p\)
for \(0<h\le h_0\).
Since \(S_h^{A_0}\) acts only on the present component,
it preserves \(\Dom(C)\).
Hence
\[
CS_h^{AC}
=
CS_h^CS_h^{A_0}
=
S_h^CCS_h^{A_0}.
\]
so
\[
\|S_h^{AC}\|_{X_\theta\to\mathcal E_p}\lesssim 1,
\qquad
\|C S_h^{AC}\|_{X_\theta\to\mathcal E_p}
=\|C S_h^{C} S_h^{A_0}\|_{X_\theta\to\mathcal E_p}\lesssim h^{-(1-\theta)}.
\]
Using the graph bound on $D(t)$ gives
\[
\|H_h(t)\|_{X_\theta\to\mathcal E_p}
=h\|D(t)S_h^{AC}\|_{X_\theta\to\mathcal E_p}
\lesssim h\|S_h^{AC}\|+h\|C S_h^{AC}\|\lesssim h^\theta.
\]
The right-preconditioned defect identity then yields the claim under Assumption~\ref{ass:OS_sectorial}
and the uniform bound on $R_h(t)$.
\end{proof}

\begin{remark}[Weaker regime: merely $\Dom(C)$--bounded $D(t)$]
\label{rem:sectorial_weaker_dom_bounded}
If one assumes only $D(t)\in\mathcal L(\Dom(C),\mathcal E_p)$ uniformly in $t$ (without graph control),
then the best general estimate for the preconditioned size is
\[
\|H_h(t)\|_{\mathcal L(X_\theta,\mathcal E_p)}\lesssim h^{1-\theta},
\]
hence the local defect bound typically degrades to
\[
\|R_h(t)-P_h(t)\|_{\mathcal L(X_\theta,\mathcal E_p)}\lesssim h^{1-\theta},
\]
uniformly in $t\in[0,T]$ (under Assumption~\ref{ass:OS_sectorial} and boundedness of $R_h(t)$).
\end{remark}
\begin{remark}
The next proposition is technical. The additional regularity assumption in 
\[
H_h(t):\Dom(C^2)\to\Dom(C)
\]
reflects the fact that the second-order defect estimate in the Proposition \ref{prop:sectorial_E_h_estimate_second_order}
requires one additional transport derivative.
\end{remark}

\begin{prop}[Second-order defect estimate on $\Dom(C^2)$ in sectorial block model]
\label{prop:sectorial_E_h_estimate_second_order}
Set
\[
S_h^{AC}=(I-hC)^{-1}(I-hA_0)^{-1},
\qquad
H_h(t)=hD(t)S_h^{AC},
\qquad
Q_h(t)=(I-H_h(t))^{-1}.
\]

Assume in addition that
\[
H_h(t):\Dom(C^2)\to\Dom(C)
\]
and
\[
\|H_h(t)\|_{\mathcal L(\Dom(C^2),\Dom(C))}
\le C h
\]
uniformly for $0<h\le h_0$:
\[
\sup_{t\in[0,T]}\sup_{0<h\le h_0}
\|Q_h(t)\|_{\mathcal L(\Dom(C))}<\infty.
\]

Then there exists $C>0$ such that
\[
\|E_h(t)x\|_{\mathcal E_p}
\le
Ch^2\|x\|_{\Dom(C^2)}.
\]
\end{prop}
\begin{proof}
Since
\[
R_h=(I-h(A_0+C+D))^{-1}
=
S_h^{AC}(I-H_h)^{-1}
=
S_h^{AC}Q_h.
\]
and
\[
P_h(t)=Q_h(t)S_h^{AC},
\]
we obtain
\[
E_h(t)
=
R_h(t)-P_h(t)
=
S_h^{AC}Q_h(t)-Q_h(t)S_h^{AC}.
\]
Hence
\[
E_h(t)
=
(S_h^{AC}-I)Q_h(t)H_h(t)
-
Q_h(t)H_h(t)(S_h^{AC}-I).
\]

Since
\[
S_h^{AC}=S_h^CS_h^{A_0},
\]
we decompose
\[
S_h^{AC}-I
=
(S_h^C-I)S_h^{A_0}
+
(S_h^{A_0}-I).
\]
Using
\[
S_h^C-I=hCS_h^C,
\qquad
S_h^{A_0}-I=hA_0S_h^{A_0},
\]
together with the boundedness of
\[
S_h^{A_0}
\]
and the standard sectorial resolvent estimate
\[
\|A_0S_h^{A_0}\|_{\mathcal L(\mathcal E_p)}
\lesssim h^{-1},
\]
we obtain
\[
\|S_h^{AC}-I\|_{\Dom(C)\to\mathcal E_p}
\lesssim h.
\]

Next we estimate $H_h(t)$ on $\Dom(C)$.
Since $D(t)$ is graph-bounded with respect to $C$,
\[
\|D(t)y\|
\lesssim
\|y\|+\|Cy\|,
\qquad y\in\Dom(C),
\]
and since
\[
S_h^{AC}x\in\Dom(C),
\qquad
CS_h^{AC}x=S_h^{AC}Cx,
\]
we obtain, for $x\in\Dom(C)$,
\[
\begin{aligned}
\|H_h(t)x\|
&=
h\|D(t)S_h^{AC}x\| \\
&\lesssim
h\bigl(
\|S_h^{AC}x\|
+
\|CS_h^{AC}x\|
\bigr) \\
&\lesssim
h\bigl(
\|x\|
+
\|Cx\|
\bigr).
\end{aligned}
\]
Thus
\[
\|H_h(t)\|_{\Dom(C)\to\mathcal E_p}
\lesssim h.
\]

Now let $x\in\Dom(C^2)$.
Since
\[
S_h^{AC}-I
=
(S_h^C-I)S_h^{A_0}
+
(S_h^{A_0}-I),
\]
and both terms map \(\Dom(C^2)\) into \(\Dom(C)\),
we conclude that
\[
(S_h^{AC}-I)x\in\Dom(C),
\qquad x\in\Dom(C^2).
\]
Assume that \(S_h^{A_0}\) preserves \(\Dom(C^2)\) uniformly in \(h\).
Using the assumed boundedness of $Q_h(t)$ on $\Dom(C)$,
together with
\[
H_h(t):\Dom(C)\to\mathcal E_p,
\]
we obtain
\[
\|Q_h(t)H_h(t)(S_h^{AC}-I)x\|
\lesssim
h^2\|x\|_{\Dom(C^2)}.
\]

By assumption,
\[
H_h(t):\Dom(C^2)\to\Dom(C)
\]
with operator norm \(O(h)\).
Hence, for \(x\in\Dom(C^2)\),
\[
H_h(t)x\in\Dom(C).
\]
Since \(Q_h(t)\) acts boundedly on \(\Dom(C)\),
we conclude that
\[
Q_h(t)H_h(t)x\in\Dom(C).
\]

Moreover, by assumption,
\[
\|H_h(t)\|_{\Dom(C^2)\to\Dom(C)}
\lesssim h.
\]

Similarly, since
\[
Q_h(t)H_h(t)x\in\Dom(C),
\]
and
\[
\|S_h^{AC}-I\|_{\Dom(C)\to\mathcal E_p}\lesssim h,
\]
we obtain
\[
\|(S_h^{AC}-I)Q_h(t)H_h(t)x\|
\lesssim
h^2\|x\|_{\Dom(C^2)}.
\]
\end{proof}
\begin{remark}[\(S_h^{A_0}\) preserves \(\Dom(C^2)\) uniformly in \(h\).]
Although we assume that \(S_h^{A_0}\) preserves \(\Dom(C^2)\) uniformly in \(h\),
since \(A_0\) acts only on the present component,
\(S_h^{A_0}\) preserves the transport regularity encoded in \(\Dom(C)\) and \(\Dom(C^2)\).
\end{remark}
\subsection{Main convergence theorem (sectorial block case)}
\label{subsec:MainC}
The following theorem shows that the sectorial block case has the same
global propagation mechanism as the non-sectorial case once the local
defect is estimated on the higher regularity scale \(Y_\theta\).
\begin{thm}[Global convergence for the sectorial block model]
\label{thm:convergence_sectorial}

Let
\[
Y_\theta:=(\Dom(C),\Dom(C^2))_{\theta,1},
\qquad 0<\theta<1.
\]
For $0<h\le h_0$ define
\[
S_h^{AC}:=(I-hC)^{-1}(I-hA_0)^{-1},
\qquad
H_h(t):=hD(t)S_h^{AC},
\]
\[
P_h(t):=(I-H_h(t))^{-1}S_h^{AC},
\qquad
R_h(t):=(I-h(A_0+C+D(t)))^{-1}.
\]

Assume the frozen local defect estimate
\[
\sup_{t\in[0,T]}
\|R_h(t)-P_h(t)\|_{\mathcal L(Y_\theta,\mathcal E_p)}
\le C h^{1+\theta}.
\]

\noindent
\textup{(i) Autonomous case.}
If $D(t)\equiv D$, assume that for every $T>0$ there exist constants
$M_{R,T}$ and $M_{P,T,\theta}$ such that
\[
\sup_{0<h\le h_0}\sup_{nh\le T}
\|R_h^n\|_{\mathcal L(\mathcal E_p)}
\le M_{R,T},
\]
and
\[
\sup_{0<h\le h_0}\sup_{nh\le T}
\|P_h^n\|_{\mathcal L(Y_\theta)}
\le M_{P,T,\theta}.
\]
Then there exists $C_{T,\theta}>0$ such that
\[
\|R_h^n-P_h^n\|_{\mathcal L(Y_\theta,\mathcal E_p)}
\le C_{T,\theta}h^\theta,
\qquad nh\le T.
\]

\noindent
\textup{(ii) Non-autonomous case.}
For $t_k=kh$, set
\[
R_k:=R_h(t_k),\qquad P_k:=P_h(t_k).
\]
Define
\[
\mathcal R_{j:k}:=R_jR_{j-1}\cdots R_k,
\qquad
\mathcal P_{j:k}:=P_jP_{j-1}\cdots P_k.
\]
Assume that
\[
\sup_{\substack{0\le k\le j\\ jh\le T}}
\|\mathcal R_{j:k}\|_{\mathcal L(\mathcal E_p)}
\le M_{R,T},
\]
and
\[
\sup_{0\le kh\le T}
\|\mathcal P_{k-1:0}\|_{\mathcal L(Y_\theta)}
\le M_{P,T,\theta}.
\]
Then there exists $C_{T,\theta}>0$ such that
\[
\|\mathcal R_{n-1:0}-\mathcal P_{n-1:0}\|_{\mathcal L(Y_\theta,\mathcal E_p)}
\le C_{T,\theta}h^\theta,
\qquad nh\le T.
\]
\end{thm}
\begin{proof}
We first prove the autonomous assertion.
Set
\[
E_h:=R_h-P_h.
\]
The telescoping identity gives
\[
R_h^n-P_h^n
=
\sum_{k=0}^{n-1}
R_h^{n-1-k}E_hP_h^k.
\]
Taking norms from $Y_\theta$ to $\mathcal E_p$, we obtain
\[
\begin{aligned}
\|R_h^n-P_h^n\|_{Y_\theta\to\mathcal E_p}
&\le
\sum_{k=0}^{n-1}
\|R_h^{n-1-k}\|_{\mathcal L(\mathcal E_p)}
\|E_h\|_{\mathcal L(Y_\theta,\mathcal E_p)}
\|P_h^k\|_{\mathcal L(Y_\theta)}.
\end{aligned}
\]
By the assumed stability bounds and the local defect estimate,
\[
\|R_h^n-P_h^n\|_{Y_\theta\to\mathcal E_p}
\le
M_{R,T}M_{P,T,\theta}
\sum_{k=0}^{n-1}
C h^{1+\theta}.
\]
Since $nh\le T$, we have $n\le T/h$. Hence
\[
\|R_h^n-P_h^n\|_{Y_\theta\to\mathcal E_p}
\le
C M_{R,T}M_{P,T,\theta} \frac{T}{h} h^{1+\theta}
\le
C_{T,\theta}h^\theta.
\]

We next prove the non-autonomous assertion.
Let
\[
E_k:=R_k-P_k.
\]
By the time-ordered telescoping identity,
\[
\mathcal R_{n-1:0}-\mathcal P_{n-1:0}
=
\sum_{k=0}^{n-1}
\mathcal R_{n-1:k+1}E_k\mathcal P_{k-1:0}.
\]
For $x\in Y_\theta$, using the stability assumptions and the uniform
frozen local defect bound, we get
\[
\begin{aligned}
\|(\mathcal R_{n-1:0}-\mathcal P_{n-1:0})x\|_{\mathcal E_p}
&\le
\sum_{k=0}^{n-1}
\|\mathcal R_{n-1:k+1}\|_{\mathcal L(\mathcal E_p)}
\|E_k\|_{\mathcal L(Y_\theta,\mathcal E_p)}
\|\mathcal P_{k-1:0}x\|_{Y_\theta}  \\
&\le
M_{R,T}\, C h^{1+\theta}\, M_{P,T,\theta}
\sum_{k=0}^{n-1}\|x\|_{Y_\theta}.
\end{aligned}
\]
Again $n\le T/h$, and therefore
\[
\|(\mathcal R_{n-1:0}-\mathcal P_{n-1:0})x\|_{\mathcal E_p}
\le
C_{T,\theta}h^\theta\|x\|_{Y_\theta}.
\]
Taking the supremum over $\|x\|_{Y_\theta}\le1$ gives
\[
\|\mathcal R_{n-1:0}-\mathcal P_{n-1:0}\|_{\mathcal L(Y_\theta,\mathcal E_p)}
\le C_{T,\theta}h^\theta.
\]
This completes the proof.
\end{proof}

\section{Numerical experiments}
\label{sec:numerical_experiments}

The purpose of the numerical experiment is not to prove the
interpolation estimate, but to illustrate the qualitative mechanism
behind the theory: higher regularity of the history variable improves
the one-step defect and prevents visible accumulation of splitting
errors over long time intervals.

The delay term was implemented through
a standard ring-buffer realization of the history variable.
At each time step, the oldest history value is discarded
and the newly computed value is appended.

\subsection{Non-sectorial model}
We consider the scalar DDE:
\[
\begin{aligned}
u'(t) &= a u(t) + b u(t + \tau), \text{ for } 0\le t\le T,\\
u(t) &= \phi(t), \text{ for }\tau\le t \le 0.
\end{aligned}
\]
We set the parameters: $a = -0.5, b = -2.0$, and $\tau =-1, T=10.0$.
First, we analyze the global error between the implicit Euler method and the Lie--Trotter splitting when the history function is sufficiently regular. In this case, $\phi\in W^{2,p}[-1, 0]$ for $p>1$.

Set the history function $\phi(t):=\cos(\pi t).$ The estimated error $E_h$: 
\[
\|E_h\|=\max_{0\le nh\le T}|u^n_R - u^n_P|
\]
is calculated in Table \ref{tab:error_on_h_cos} for $h = 2^{-6}, 2^{-7}, 2^{-8}, 2^{-9}$, respectively.

For the smooth history function \(\phi(t)=\cos(\pi t)\),
the observed convergence rate is approximately first order.
This is consistent with the theoretical mechanism:
although the one-step defect is of second order on sufficiently
regular data, the global error accumulates over \(O(h^{-1})\)
time steps, leading to an overall \(O(h)\) rate.
\begin{table}[htbp]
\begin{center}
\caption{The error $E_h$ on $h=2^{-k}$, $k=6, 7, 8, 9$.}
  \label{tab:error_on_h_cos} 
\begin{tabular}{ccc}
  \hline
  $k$ & $2^{-k}$ & Error \\
  \hline
  6 &0.015625 &0.4300819189508977 \\
  7 &0.0078125 &0.21770688077278366 \\
  8 &0.00390625 &0.10948946489554923 \\
  9 &0.001953125 &0.05489986745744313 \\
  \hline
  estimated theta & & 0.9900805672556413\\
  \hline
\end{tabular}
\end{center}
\end{table}

\begin{figure}[H]
  \centering
  \begin{tabular}{cc}
    \begin{minipage}[t]{0.48\textwidth}
      \centering
      \includegraphics[width=\linewidth]{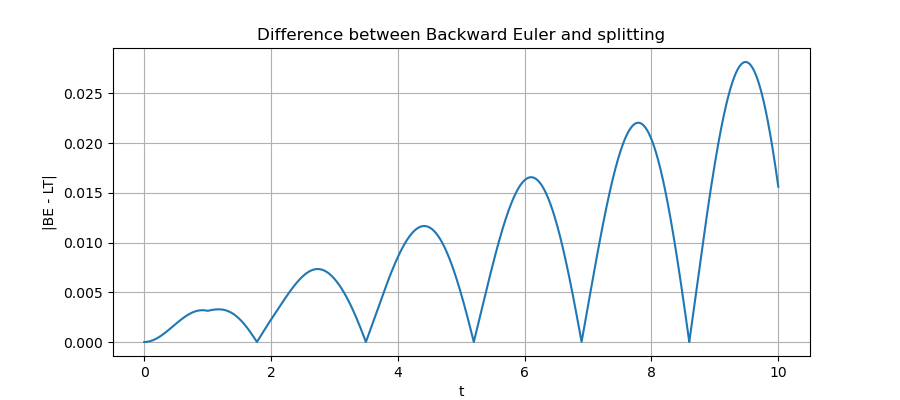}
      \caption{The error between IE and LT, history = $\cos(\pi t)$.}
      \label{Fig_error_IE_LT_hist_cosine}
    \end{minipage} &
    \begin{minipage}[t]{0.48\textwidth}
      \centering
      \includegraphics[width=\linewidth]{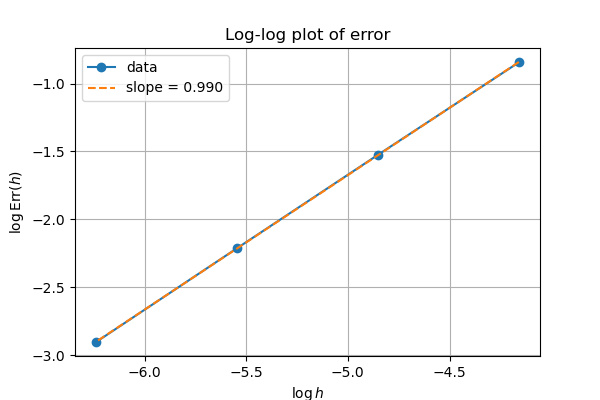}
      \caption{Log-Log Plot of the error.}
      \label{Fig_loglog_error_IE_LT_hist_cosine}
    \end{minipage}
  \end{tabular}
\end{figure}

Fig. \ref{Fig_error_IE_LT_hist_cosine} illustrates the error between Implicit Euler and Lie--Trotter for $0\le t\le T=10.0$, step-size = $h=0.001$, with history fuction: $\cos(\pi t)$, for $-1.0\le t\le 0$. 
Fig. \ref{Fig_loglog_error_IE_LT_hist_cosine} shows $\log \mathrm{Error}(h)$ vs $\log (h)$. 
\[
\log \mathrm{Err}(h)
\text{ vs }
\log h.
\]

If
\[
\mathrm{Err}(h)\approx C h^\theta
\]
then, taking the logarithm of both sides, we obtain
\[
\log \mathrm{Err}(h)\approx\log C + \theta \log h;
\]
the slope of the log-log plot is $\theta$.

The following two examples analyze a case where, if an initial history function of low regularity is in $\Dom(C)$, then even if it is not in $\Dom(C^2)$,
\[
\|E_h\|=\max_{0\le nh\le T}|u^n_R - u^n_P|\sim O(h)
\]
holds.

Let the history function be defined as $\phi(t) = |t + 1/2|^\beta$.
\begin{itemize}
 \item[(i)] $\beta = 0.5$
\begin{table}[htbp]
\begin{center}
\caption{The error $E_h$ with singular history $|1+t/2|^{0.5}$ for $h=2^{-k}$, $k=6, 7, 8, 9$.}
\label{tab:error_on_h_abs_half} 
\begin{tabular}{ccc}
  \hline
  $k$ & $2^{-k}$ & Error \\
  \hline
  6 &0.015625 &0.134930330118706 \\
  7 &0.0078125 &0.06826565416150679 \\
  8 &0.00390625 &0.034325345432295234 \\
  9 &0.001953125 &0.017209468778896814 \\
  \hline
  estimated theta & & 0.9904706299032435 \\
  \hline
\end{tabular}
\end{center}
\end{table}
\begin{figure}[H]
  \centering
  \begin{tabular}{cc}
    \begin{minipage}[t]{0.48\textwidth}
      \centering
      \includegraphics[width=\linewidth]{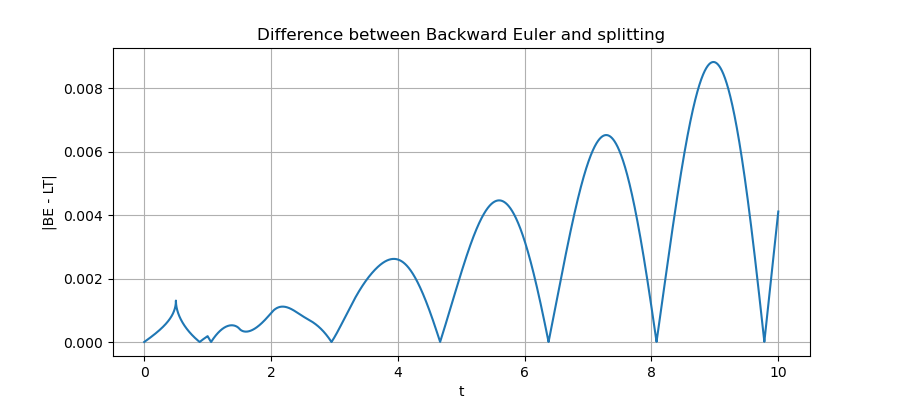}
      \caption{The error between IE and LT, history = $|t + 1/2|^{0.5}$.}
      \label{Fig_error_IE_LT_hist_signular_beta1half}
    \end{minipage} &
    \begin{minipage}[t]{0.48\textwidth}
      \centering
      \includegraphics[width=\linewidth]{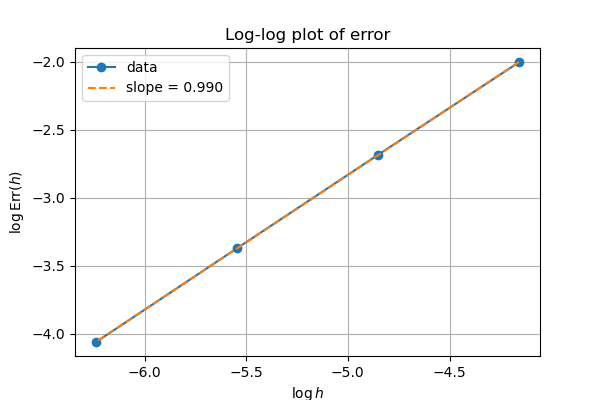}
      \caption{$\log(\mathrm{Err}(h))$ vs $\log (h)$ with singular history $=|t+1/2|^{0.5}$.}
      \label{Fig_loglog_error_IE_LT_hist_singular_beta1half}
    \end{minipage}
  \end{tabular}
\end{figure}
\item[(ii)] $\beta = 1.5$
\begin{table}[htbp]
\caption{The error $E_h$ with singular history $|t+1/2|^{1.5}$ for $h=2^{-k}$, $k=6, 7, 8, 9$.}
\label{tab:error_on_h_abs_3over2} 
\begin{center}
\begin{tabular}{ccc}
  \hline
  $k$ & $2^{-k}$ & Error \\
  \hline
  6 &0.015625 &0.05752373921901938 \\
  7 &0.0078125 &0.029065350705435866 \\
  8 &0.00390625 &0.014605914279448504 \\
  9 &0.001953125 &0.0073209815172172454 \\
  \hline
  estimated theta && 0.9914893060805645 \\
  \hline
\end{tabular}
\end{center}
\end{table}

\begin{figure}[H]
  \centering
  \begin{tabular}{cc}
    \begin{minipage}[t]{0.48\textwidth}
      \centering
      \includegraphics[width=\linewidth]{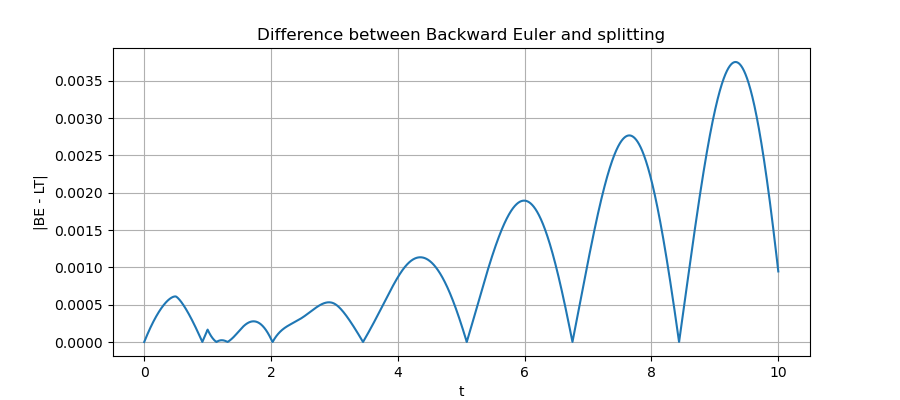}
      \caption{The error between IE and LT, history = $|t + 1/2|^{1.5}$.}
      \label{Fig_error_IE_LT_hist_signular_beta3over2}
    \end{minipage} &
    \begin{minipage}[t]{0.48\textwidth}
      \centering
      \includegraphics[width=\linewidth]{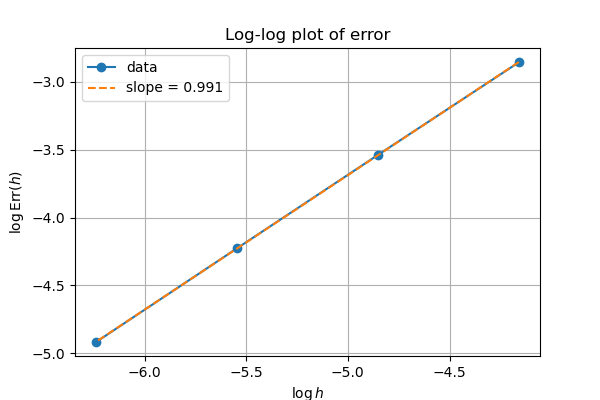}
      \caption{$\log(\mathrm{Err}(h))$ vs $\log(h)$ with singular history $=|t+1/2|^{1.5}$.}
      \label{Fig_loglog_error_IE_LT_hist_singular_beta3over2}
    \end{minipage}
  \end{tabular}
\end{figure}
\end{itemize}
Tables \ref{tab:error_on_h_abs_half} and \ref{tab:error_on_h_abs_3over2} present the error $E_h$ between Implicit Euler and Lie--Trotter, with history functions $|t+1/2|^{0.5}$ and $|t + 1/2|^{1.5}$, respectively, for $h=2^{-k}$, $k=6, 7, 8, 9$.

Figs \ref{Fig_error_IE_LT_hist_signular_beta1half} and \ref{Fig_error_IE_LT_hist_signular_beta3over2} show the errors between Implicit Euler and Lie-Trotter for $\beta =1/2$ and $\beta =3/2$, respectively.
Figs \ref{Fig_loglog_error_IE_LT_hist_singular_beta1half} and \ref{Fig_loglog_error_IE_LT_hist_singular_beta3over2} show $\log (\mathrm{Err}(h))$ vs $\log (h)$ for 
$\beta =1/2$ and $\beta =3/2$, respectively.

\subsection{Numerical observations}

The numerical experiments reveal a characteristic oscillatory structure in the splitting error. This behavior appears to be closely related to the delayed feedback mechanism inherent in the DDE. Since the difference between the implicit Euler scheme and the Lie--Trotter splitting scheme is generated by the one-step shift
\[
u^{n+1-m}-u^{n-m},
\]
the local defect may be interpreted as a discrete time derivative of the delayed term. Consequently, the resulting error is repeatedly reintroduced through the delay loop.

A notable feature of the error curves is the occurrence of repeated troughs at times close to integer multiples of the delay. Numerically, this suggests a partial resynchronization of the delayed states. At the same time, the slowly increasing envelope indicates the cumulative effect of local splitting defects over successive delay cycles.

Another observation is that the overall oscillatory pattern is remarkably robust with respect to the regularity of the history function. Both smooth and singular histories exhibit essentially the same delay-induced structure. The principal influence of the history regularity appears in the amplitude and local shape of the error rather than in its qualitative long-time behavior.

These observations are consistent with the theoretical viewpoint developed in this paper: the dominant features of the splitting error are determined by the interaction between local defects, delayed feedback, and discrete stability of the associated propagators.

The singular histories provide an interesting comparison with the smooth case.
For
\[
\phi(t)=|t+1/2|^{3/2},
\]
the history belongs to $W^{1,p}[-1, 0]$ for $p\ge 2$ but not to $W^{2,p}[-1,0]$, which is precisely the type of regularity considered in the interpolation framework developed in this paper. Nevertheless, the observed convergence rate remains close to first order.

For the more singular history
\[
\phi(t)=|t+1/2|^{1/2},
\]
the same qualitative behavior is observed numerically, although this example lies partially outside the regularity assumptions used in the theoretical analysis, since this $\phi(t)$ does not even belong to $W^{1,p}[-1, 0]$, for $p\ge 2$. 

In both cases the error decays approximately as $O(h)$, indicating that the discrete propagators remain stable even when the history variable possesses limited smoothness.

\subsection{Sectorial block model}
As a second experiment, we consider a sectorial block example
given by a delayed reaction--diffusion equation. The purpose is not
to provide a full numerical study of delay PDEs, but to verify that
the same discrete-resolvent mechanism is observed when the transport
operator is coupled with a sectorial spatial operator.

\paragraph{Model and discretization.}
We consider the reaction--diffusion delay equation on $(0,1)$ with Dirichlet boundary conditions:
\begin{equation}\label{eq:num_reaction_diffusion}
\begin{aligned}
u_t(t,x) &= \kappa u_{xx}(t,x)+ a u(t,x)+ b\,u(t+\tau,x),
\qquad t\in[0,T],\ x\in(0,1),\\
u(t,0)&=u(t,1)=0,\qquad t\ge0,\\
u(t,x)&=\phi(t,x)=\cos(\pi t)\sin(\pi x),\qquad t\in[\tau,0].
\end{aligned}
\end{equation}
The parameters $a, b, T, \tau$ are the same as in the non-sectorial model above. We measure the error between IE vs LT in $L^2(0, 1)$.

The Implicit Euler step is
\[
(I-h(\kappa \Delta_h+aI))u_R^{n+1} = u_R^n + h b\,u_R^{n+1-m}.
\]
The Lie--Trotter splitting step is
\[
(I-h(\kappa \Delta_h+aI))u_P^{n+1} = u_P^n + hb u_P^{n-m}.
\]
The difference between them is also
\[
u^{n+1-m} \quad\text{vs}\quad u^{n-m}.
\]

We computed the error: $E_h$ in $L^2(0,1)$.
The result is as follows:
\begin{table}
\begin{center}
\caption{Reaction-diffurion $L^2(0,1)$ error between Implicit Euler and Lie-Trotter, for $h=2^{-k}$, $k=6, 7, 8, 9$.}
\label{tab:error_IE_vs_LT_react_diff}
\begin{tabular}{ccc}
  \hline
  $k$ & $2^{-k}$ & Error \\
  \hline
  6 &0.015625 &0.05752373921901938 \\
  7 &0.0078125 &0.029065350705435866 \\
  8 &0.00390625 &0.014605914279448504 \\
  9 &0.001953125 &0.0073209815172172454 \\
  \hline
  estimated theta & & 0.9914893060805645 \\
  \hline
\end{tabular}
\end{center}

\end{table}
\begin{figure}[H]
  \centering
    \includegraphics[width=\linewidth]{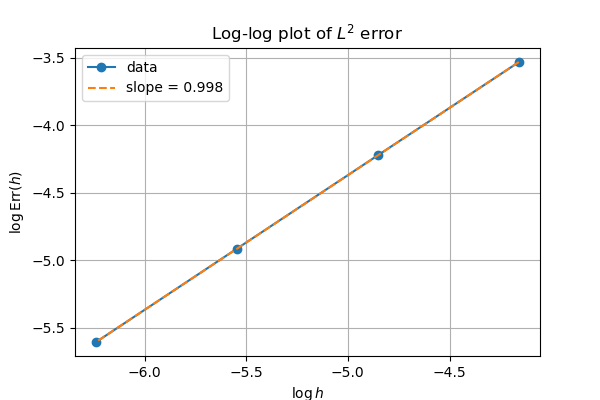}
      \caption{The error in $L^2[0,1]$ between IE and LT, history = $\cos(\pi t)\sin(\pi x)$.}
      \label{Fig_loglog_error_react_diff_L2}
\end{figure}
We obtain the global error $\sim O(h)$ as $h\to 0$: 
\[
\mathrm{Err}(h)
=
\max_{0\le nh\le T}
\|u_R^n-u_P^n\|_{L^2(0,1)}\sim O(h) \text{ as }h\to 0.
\]

The sectorial reaction–diffusion example exhibits essentially the same convergence behavior as the scalar delay equation. The observed rate is again close to first order, suggesting that the introduction of a sectorial spatial operator does not alter the fundamental local-to-global mechanism of the discrete-resolvent framework.

These results support the theoretical viewpoint that the dominant contribution to the splitting error originates from the delay interaction and its propagation through the discrete resolvent structure.
\subsection{Conclusion of the numerical experiments.}
\paragraph{Conclusion of the numerical experiments.}

The numerical experiments support the qualitative picture developed in the theoretical part of the paper. For both smooth and singular histories, the difference between the implicit Euler scheme and the Lie--Trotter splitting scheme decreases approximately linearly with respect to the stepsize.

In particular, no visible loss of convergence is observed for histories of limited regularity. Moreover, the same behavior appears in the delayed reaction--diffusion model, indicating that the local-to-global mechanism identified in the discrete-resolvent framework persists in both non-sectorial and sectorial settings.

While the numerical results do not constitute a proof of the interpolation estimates, they provide evidence that the discrete stability mechanism developed in this paper captures the dominant features of the splitting error. The experiments also suggest that the theoretical bounds may not be sharp, since the observed convergence rate remains close to first order even in situations where only weaker regularity is available.

\section{Conclusion}
\label{sec:conclusion}
\subsection{Summary of the structural results}

This paper develops a discrete resolvent framework for the analysis of
implicit Euler and Lie--Trotter splitting schemes for delay differential equations.

The analysis is carried out entirely at the level of resolvent-type propagators of delay equations.

A central observation is that the transport resolvent
\[
S_h^C=(I-hC)^{-1}
\]
naturally induces a hierarchy of fractional interpolation spaces.
On the basic interpolation scale
\[
X_\theta=(\mathcal E_p,\Dom(C))_{\theta,1},
\]
the local defect satisfies
\[
\|R_h(t)-P_h(t)\|_{\mathcal L(X_\theta,\mathcal E_p)}
\lesssim h^\theta
\]
in the non-sectorial case.

By introducing the higher regularity scale
\[
Y_\theta=(\Dom(C),\Dom(C^2))_{\theta,1},
\]
the local defect improves to
\[
\|R_h(t)-P_h(t)\|_{\mathcal L(Y_\theta,\mathcal E_p)}
\lesssim h^{1+\theta}.
\]

The same mechanism extends to sectorial block models, where the transport
and sectorial components can be separated at the resolvent level through
an exact factorization of the $(A_0+C)$ resolvent.

\subsection{Local-to-global mechanism}

Local defect estimates alone do not automatically imply
global convergence.
Under mere power-boundedness assumptions,
the telescoping expansion produces a factor of order
\[
n\sim h^{-1},
\]
which may destroy convergence when only
\[
\|E_h\|\lesssim h^\theta
\]
is available.

The numerical experiments indicate that the observed global
difference between implicit Euler and Lie--Trotter splitting
decays approximately with first order in all tested examples,
including histories with limited regularity and a delayed
reaction--diffusion model.

Although these experiments do not constitute proof beyond
the abstract assumptions of the theory, they suggest that the
discrete-resolvent mechanism remains effective in a wider class
of examples than those covered by the convergence theorems.

\subsection{Main contribution of this paper}

The main contribution of the paper is therefore not a new splitting
scheme itself, but a discrete operator-theoretic framework that
explains how local resolvent defects and discrete stability combine
to produce global convergence in both autonomous and non-autonomous
delay equations.

\subsection{Outlook}

The present framework obtains global convergence by combining higher-order local defect estimates on $Y_\theta$ with power-bounded stability assumptions. An interesting direction for future research is to investigate whether stronger discrete stability properties, such as Ritt-type estimates, bounded-variation conditions, Abel summation techniques, or Kreiss-type bounds, allow global convergence to be established directly from the $X_\theta$-level defect estimates, without introducing the higher-regularity space $Y_\theta$.

\appendix

\appendix
\section{Transport resolvent estimates}
\label{app:transport-resolvent}

Throughout this appendix let $\tau<0$ and $1<p<\infty$, and set
\[
\mathcal E_p := X \oplus L^p([\tau,0];X),
\qquad
\|(u,\rho)\|_{\mathcal E_p} := \|u\|_X + \|\rho\|_{L^p([\tau,0];X)}.
\]
Let $C$ be the transport operator
\[
C(u,\rho) := (0,\partial_\sigma \rho),
\qquad
\Dom(C) := \{(u,\rho)\in X\oplus W^{1,p}([\tau,0];X):\rho(0)=u\}.
\]
For $h>0$ we write
\[
S_h^C := (I-hC)^{-1}\in\mathcal L(\mathcal E_p).
\]

\subsection{Explicit formula and basic bounds}
\label{app:explicit-formula}

\begin{prop}[Explicit form of $S_h^C$]
\label{prop:ShC-explicit}
Let $(f,g)\in\mathcal E_p$. Then $S_h^C(f,g)=(u,\rho)$ is given by
\[
u=f,
\qquad
\rho(\sigma)
= e^{\sigma/h} f + \frac{1}{h}\int_{\sigma}^{0} e^{(\sigma-\eta)/h}\, g(\eta)\, d\eta,
\quad \sigma\in[\tau,0].
\]
In particular, $S_h^C$ maps $\mathcal E_p$ into $\Dom(C)$.
\end{prop}

\begin{proof}
Solve $(I-hC)(u,\rho)=(f,g)$, that is,
\[
u=f,\qquad \rho-h\rho'=g,\qquad \rho(0)=u=f.
\]
Multiplying the differential equation by $e^{-\sigma/h}$ and integrating from
$\sigma$ to $0$ yields the stated formula.
\end{proof}

\subsection{Sharp trace estimate}
\label{app:trace-estimate}

\begin{prop}[Trace bound]
\label{app:prop:ShC-trace}
There exists $C>0$ (depending only on $p$ and $\tau$) such that for all $h\in(0,1]$,
\[
\|\Phi_\tau S_h^C(f,g)\|_X
\le \|f\|_X + C\, h^{-1/p}\,\|g\|_{L^p([\tau,0];X)}.
\]
\end{prop}

\begin{proof}
Let $(u,\rho)=S_h^C(f,g)$. By Proposition~\ref{prop:ShC-explicit} at $\sigma=\tau$,
\[
\rho(\tau)=e^{\tau/h}f+\frac{1}{h}\int_{\tau}^{0}e^{(\tau-\eta)/h}g(\eta)\,d\eta.
\]
Since $\tau<0$, we have $0<e^{\tau/h}\le 1$, hence $\|e^{\tau/h}f\|\le \|f\|$.
For the integral term, apply H\"older's inequality:
\[
\Bigl\|\frac{1}{h}\int_{\tau}^{0}e^{(\tau-\eta)/h}g(\eta)\,d\eta\Bigr\|
\le \Bigl\|\frac{1}{h}e^{(\tau-\cdot)/h}\Bigr\|_{L^{p'}([\tau,0])}\,
\|g\|_{L^p([\tau,0];X)}.
\]
A direct computation gives
\[
\Bigl\|\frac{1}{h}e^{(\tau-\cdot)/h}\Bigr\|_{L^{p'}([\tau,0])}^{p'}
= \int_{\tau}^{0}\Bigl(\frac{1}{h}e^{(\tau-\eta)/h}\Bigr)^{p'}\,d\eta
= h^{1-p'}\int_{0}^{-\tau/h}e^{-p's}\,ds
\le \frac{1}{p'}\,h^{1-p'}.
\]
Therefore
\[
\Bigl\|\frac{1}{h}e^{(\tau-\cdot)/h}\Bigr\|_{L^{p'}([\tau,0])}
\le (p')^{-1/p'}\,h^{-1/p},
\]
and the claim follows.
\end{proof}

\subsection{Fractional lifting via interpolation}
\label{app:fractional-estimate}

Let
\[
X_\theta := (\mathcal E_p,\Dom(C))_{\theta,1},\qquad 0<\theta<1.
\]

\begin{prop}[Fractional lifting estimate]
\label{prop:ShC-fractional}
For every $0<\theta<1$ there exists $C_\theta>0$ such that for all $h\in(0,1]$,
\[
\|S_h^C\|_{\mathcal L(X_\theta,\Dom(C))}\le C_\theta\,h^{-\theta}.
\]
\end{prop}

\begin{proof}
Endow $\Dom(C)$ with the graph norm
$\|x\|_{\Dom(C)}:=\|x\|_{\mathcal E_p}+\|Cx\|_{\mathcal E_p}$.
By \eqref{eq:resolvent_inequality}
and
\eqref{eq:C_resolvent_of_C_inequality},
\[
\|S_h^C\|_{\mathcal L(\mathcal E_p,\mathcal E_p)}
\le M_0,
\qquad
\|S_h^C\|_{\mathcal L(\mathcal E_p,\Dom(C))}
\lesssim h^{-1}.
\]
Real interpolation between the bounds $\mathcal E_p\to\mathcal E_p$ and
$\mathcal E_p\to\Dom(C)$ yields
\[
\|S_h^C\|_{\mathcal L((\mathcal E_p,\Dom(C))_{\theta,1},\,\Dom(C))}
\lesssim h^{-\theta},
\]
which is the claim.
\end{proof}
\begin{lemma}[Estimate for $C S_h^C$ on fractional spaces]
\label{lem:CS_hC_fractional}

Let $0<\theta<1$ and
$X_\theta=(\mathcal E_p,\Dom(C))_{\theta,1}$.
Then there exists $C_\theta>0$ such that for all $0<h\le 1$,
\[
\|C S_h^C\|_{\mathcal L(X_\theta,\mathcal E_p)}
\le C_\theta\, h^{-(1-\theta)}.
\]

\end{lemma}

\begin{proof}
We interpolate between the two endpoint estimates.

\medskip
\noindent
\textbf{(i) On $\mathcal E_p$.}
Since
\[
C S_h^C = \frac{1}{h}(S_h^C - I),
\]
and $\|S_h^C\|_{\mathcal L(\mathcal E_p)}\le M_0$,
we obtain
\[
\|C S_h^C\|_{\mathcal L(\mathcal E_p)}
\le \frac{2}{h}.
\]

\medskip
\noindent
\textbf{(ii) On $\Dom(C)$.}
For $x\in\Dom(C)$,
the commutation relation $C S_h^C x = S_h^C Cx$
implies
\[
\|C S_h^C x\|
\le \|S_h^C\|\,\|Cx\|
\le \|Cx\|.
\]
Hence
\[
\|C S_h^C\|_{\mathcal L(\Dom(C),\mathcal E_p)} \le 1.
\]

\medskip
\noindent
Interpolating between
$\mathcal E_p$ and $\Dom(C)$
with real interpolation
$X_\theta=(\mathcal E_p,\Dom(C))_{\theta,1}$
yields
\[
\|C S_h^C\|_{\mathcal L(X_\theta,\mathcal E_p)}
\lesssim
(h^{-1})^{1-\theta}\cdot 1^{\theta}
=
h^{-(1-\theta)}.
\]
\end{proof}



\section*{Acknowledgments}
I would like to express my sincere gratitude to my advisor, Professor Toru Ohira of the Graduate School of Mathematics at Nagoya University, for his invaluable guidance, encouragement, and insightful discussions throughout this research. I would also like to thank Professor Toshiaki Hishida of the Graduate School of Mathematics at Nagoya University for enthusiastically answering my questions and for his invaluable advice on generating semigroups and evolution families while preparing this paper. 

\bibliographystyle{plain} 
\bibliography{discrete_resolvent_DDE}
\end{document}